\documentclass{amsart}
\usepackage{amsthm,amsmath,amssymb,amsfonts, amscd,graphics, latexsym,
  enumerate, stmaryrd,xspace,verbatim,epic,eepic}
\usepackage[all]{xypic}

\SelectTips{cm}{}

\newcounter{endnotecounter}
\newcounter{endendnotecounter}
\newcommand{\note}[1]%
{\refstepcounter{endnotecounter}%
\label{Endnotes:\theendnotecounter}%
\makebox[0pt]{\raisebox{1ex}[0in][0in]{\tiny\bf{*}\theendnotecounter{*}}}%
\marginpar%
 [\hfill\rm\theendnotecounter$\rightarrow$]{$\leftarrow$\rm\theendnotecounter}%
\AtEndDocument{\stepcounter{endendnotecounter}{%
 \begin{trivlist} \item[\theendendnotecounter\  (page
 \pageref{Endnotes:\theendendnotecounter}).] #1
\end{trivlist}%
}}}

 \newlength{\baseunit}               
 \newcount{\numlines}                
{
   \newtheorem{theorem}[subsubsection]{Theorem}
   \newtheorem{proposition}[subsubsection]{Proposition}
   \newtheorem{lemma}[subsubsection]{Lemma}

   \newtheorem{corollary}[subsubsection]{Corollary}

}
{\theoremstyle{definition}

   \newtheorem{definition}[subsubsection]{Definition}
   \newtheorem{remark}[subsubsection]{Remark}
   
}
\newcommand{\bbS}{{\mathbb{S}}}

\newcommand{\CC}{{\mathbb{C}}}
\newcommand{\QQ}{{\mathbb{Q}}}

\newcommand{\PP}{{\mathbb{P}}}
\newcommand{\ZZ}{{\mathbb{Z}}}

\newcommand{\GG}{{\mathbb{G}}}

\newcommand{\bbA}{{\mathbb{A}}}
\renewcommand{\AA}{{\mathbb{A}}}

\newcommand{\bG}{{\mathbf{G}}}

\newcommand{\bL}{{\mathbf{L}}}

\newcommand{\bR}{{\mathbf{R}}}

\newcommand{\bmu}{{\boldsymbol{\mu}}}
\newcommand{\bone}{{\boldsymbol{1}}}

\newcommand{\fr}{{\mathfrak{r}}}

\newcommand{\cA}{{\mathcal A}}
\newcommand{\cB}{{\mathcal B}}
\newcommand{\cC}{{\mathcal C}}
\newcommand{\cE}{{\mathcal E}}
\newcommand{\cF}{{\mathcal F}}
\newcommand{\cG}{{\mathcal G}}
\renewcommand{\cH}{{\mathcal H}}
\newcommand{\cI}{{\mathcal I}}
\newcommand{\cK}{{\mathcal K}}
\renewcommand{\cL}{{\mathcal L}}
\newcommand{\cM}{{\mathcal M}}

\newcommand{\cO}{{\mathcal O}}

\newcommand{\cT}{{\mathcal T}}
\newcommand{\cU}{{\mathcal U}}
\newcommand{\cV}{{\mathcal V}}

\newcommand{\cX}{{\mathcal X}}
\newcommand{\cY}{{\mathcal Y}}

\newcommand{\fB}{{\mathfrak B}}

\newcommand{\fC}{{\mathfrak C}}
\newcommand{\fD}{{\mathfrak D}}
\newcommand{\fM}{{\mathfrak M}}
\newcommand{\fP}{{\mathfrak P}}
\newcommand{\fQ}{{\mathfrak Q}}
\newcommand{\tw}{{\operatorname{tw}}}
\newcommand{\twD}{{\fD^{\operatorname{tw}}}}
\newcommand{\twM}{{\fM^{\operatorname{tw}}}}

\newcommand{\vir}{{\text{vir}}}

\newcommand{\Spec}{\operatorname{Spec}}
\newcommand{\Spf}{\operatorname{Spf}}

\newcommand{\Isom}{\operatorname{Isom}}

\newcommand{\Sing}{\operatorname{Sing}}

\newcommand{\Pic}{{\operatorname{\mathbf{Pic}}}}
\newcommand{\Hom}{{\operatorname{Hom}}}
\newcommand{\cHom}{{{\cH}om}}

\newcommand{\Aut}{{\operatorname{Aut}}}
\newcommand{\cAut}{{{\cA}ut}}

\newcommand{\codim}{\operatorname{codim}}

\newcommand{\lrar}{\longrightarrow}

\newcommand{\ocM}{\overline{{\mathcal M}}}

\newcommand{\proj}{{\mathbb P}}

\newcommand{\double}{\genfrac..{0pt}1
{\raise -3pt\hbox{$\scriptstyle\longrightarrow$}}{\raise 5pt\hbox
{$\scriptstyle\longrightarrow$}}}
\newcommand{\setmin}{\,\protect%
\begin{picture}(8,3.5)\qbezier(1,3.5)(4,2.)(7,.5)\end{picture}\,}

\renewcommand{\setminus}{\setmin}

\newcommand{\hide}[1]{}

\renewcommand\H{\operatorname{H}}
\newcommand\QH{\operatorname{QH}}
\newcommand\QA{\operatorname{QA}}

\def\eqdef{\mathrel{\mathop=\limits^{\rm def}}}

\newcommand\noqed{\renewcommand\qed{}}

\def\tototi{\mathbin{\mathop{\otimes}\limits^{\raise-1pt\hbox
{$\scriptscriptstyle {\rm L}$}}}}

\def\indlim{\mathop{\vrule width0pt height7pt depth
4pt\smash{\lim\limits_{\raise 1pt\hbox to 14.5pt
{\rightarrowfill}}}}}
\def\projlim{\mathop{\vrule width0pt height7pt depth
4pt\smash{\lim\limits_{\raise 1pt\hbox to 14.5pt
{\leftarrowfill}}}}}

\newcommand\radice[2]{\sqrt[\uproot{2}#1]{#2}}
\renewcommand\th{^{\text{th}}}
\newcommand\gm{\GG_\mathrm{m}}

\newcommand\displaceamount{2.8pt}

\newcommand{\doubledown}{\ar@<\displaceamount>[d]\ar@<-\displaceamount>[d]}

\newcommand{\doubleup}{\ar@<\displaceamount>[u]\ar@<-\displaceamount>[u]}

\newcommand{\doubleright}{\ar@<\displaceamount>[r]\ar@<-\displaceamount>[r]}

\newcommand{\iX}[1][]{{\cI_{\bmu_{#1}}(\cX)}}
\newcommand{\riX}[1][]{{\overline{\cI}_{\bmu_{#1}}(\cX)}}

\newcommand{\tworiX}[1][]
   {{\overline{\cI}^{(2)}_{\bmu_{#1}}(\cX)}}
\newcommand{\riXs}{{\riX^2}}

\newcommand{\checkbullet}{{\check\bullet}}

\newcommand{\thickslash}{\mathbin{\!\!\pmb{\fatslash}}}
\newcommand{\slashprime}{\thickslash_{1}\!}
\newcommand{\slashsecond}{\thickslash_{2}\!}
\newcommand{\age}{\mathrm{age}}
\newcommand{\cyc}{\mathrm{cyc}}
\newcommand{\PD}{\mathrm{PD}}
\newcommand{\tor}{\operatorname{Tor}}
\newcommand{\rk}{\operatorname{rk}}

\def\into{\hookrightarrow}

\newcommand{\stringy}{stringy\xspace}
\newcommand{\stChow}[2][*]{{A^{#1}_{\mathrm{st}}(#2)_{\QQ} }}
\newcommand{\oh}{{\mathcal O}}

\newcommand\ddash{\nobreakdash--\hspace{0pt}}

\allowdisplaybreaks

\begin{document}

\title[Gromov--Witten theory of stacks]{Gromov--Witten theory\\of Deligne--Mumford stacks}

\author[D. Abramovich]{Dan Abramovich}
\thanks{Research of D.A. partially supported by NSF grant DMS-0335501}
\address{Department of Mathematics, Box 1917, Brown University,
Providence, RI, 02912, U.S.A}
\email{abrmovic@math.brown.edu}
\author[T. Graber]{Tom Graber}
\thanks{Research of T.G. partially supported by NSF grant DMS-0301179
and an Afred P. Sloan research fellowship.}
\address{Department of Mathematics\\
California Institute of Technology\\
Pasadena, CA \\
 U.S.A.} \email{graber@caltech.edu}
\author[A. Vistoli]{Angelo Vistoli}
\thanks{Research of A.V. partially supported by the PRIN Project ``Geometria
sulle variet\`a algebriche'', financed by MIUR}
\address{Scuola Normale Superiore\\ Piazza dei Cavalieri 7\\
56126 Pisa\\ Italy}
\email{angelo.vistoli@sns.it}



\maketitle

\setcounter{tocdepth}{1}
\tableofcontents
\section{Introduction}
Gromov\ddash Witten theory of orbifolds was introduced in the symplectic
setting in \cite{Chen-Ruan}. In \cite{AGV} we adapted the theory to
algebraic geometry, using Chow rings and the language of stacks. The
latter work is, to a large extent, a research announcement, as
detailed  proofs were not given.   The main purpose of this paper is
to complete that work and lay the algebro-geometric theory on a
sounder footing.

It appears from the emerging literature that the language of stacks
-- whether algebraic or differential -- is imperative in making
serious computations in orbifold Gromov\ddash Witten theory (see, e.g.
\cite{Cadman2}, \cite{Tseng}). It can therefore be hoped that some of
the material here should be useful in the symplectic setting as
well.

The fact that a few years have passed since the release of our paper
\cite{AGV} may be one cause for a regretful slow development of
applications. On the other hand, we believe recent developments
have enabled us to set the foundations in a much better way than was
possible at the time of  the paper \cite{AGV}. Most important among
these are Olsson's papers \cite{Olsson-Hom} and
\cite{Olsson-twisted}.

\subsection{Twisted stable maps}
Algebro-geometric treatments of Gromov\ddash Witten theory of a smooth
projective variety $X$ rely on Kontsevich's moduli stacks
$\ocM_{g,n}(X, \beta)$, parametrizing $n$-pointed stable maps from
curves of genus $g$ to $X$ with image class $\beta\in \H_2(X, \ZZ)$,
see  \cite{BM}. When one replaces the manifold $X$ with an
orbifold, one needs to replace the curves in the stable maps by
orbifold curves - this is a phenomenon discussed in detail elsewhere
(see \cite{AV-families}). The stack of stable maps is thus replaced
by the stack of \emph{twisted stable maps}, denoted $\cK_{g,n}(\cX,
\beta)$ in \cite{AV}, where it was constructed.

The proof of the main theorem in \cite{AV} is not ideal
as it relies on ad-hoc arguments and requires verifying Artin's
axioms one by one. An alternative approach which is much more
conceptual is given in Olsson's papers  \cite{Olsson-Hom} and
\cite{Olsson-twisted}. Artin's axioms still need to be verified, but
in a more general and cleaner situation. In Olsson's paper
\cite{Olsson-twisted} one also finds a direct construction of the
stack of twisted pre-stable curves, which is an important tool in
the theory developed here.  See Section \ref{Sec:twisted-curves} (and
Appendix \ref{roots}) for a quick review.

An analogous space of \emph{orbifold stable maps} was constructed by
W. Chen and Y. Ruan in \cite{Chen-Ruan} using very different
methods. In spirit the two constructions describe the same thing, and one expects that the resulting Gromov--Witten numbers are identical.

For quotient stacks by a finite group $\cX= [V/G]$, Jarvis, Kaufmann
and Kimura considered in \cite{JKK} maps of pointed admissible
$G$-covers to  $V$.  This was revisited in \cite{AGOT}. Lev
Borisov recently discovered in his mail archives two letters from
Kontsevich, dated from July 1996, where Gromov--Witten theory of  a
quotient stack by a finite group is outlined precisely using
admissible $G$-covers. See \cite{cimenotes} for a reproduced text.

When the target stack $\cX$ is smooth, the stack $\cK_{g,n}(\cX,
\beta)$ admits a perfect obstruction theory in the sense of
\cite{BF} and so one gets a virtual fundamental class
$[\cK_{g,n}(\cX, \beta)]^\vir$, which then gives rise to a
Gromov-Witten theory.

\subsection{Gromov--Witten classes - manifold case}
We follow a formalism for Gromov--Witten theory suited for Chow
rings developed in \cite{GP} (a paper longer
overdue than this one). A similar formalism was given in \cite{EK}, Section~3. Consider classes $\gamma_1,\ldots,\gamma_n
\in A^*(X)$. Define the associated Gromov--Witten classes to be
$$\langle \gamma_1,\ldots,\gamma_n ,*\rangle^X_{g,\beta} \ \ = \\  e_{n+1\ *}(e_1^*\gamma_1 \cup\ldots \cup e_n^*\gamma_n\cap[\ocM_{g,n+1}(X,\beta)]^\vir)\in A^*(X).$$ An important part of Gromov--Witten theory is concerned with relations these classes satisfy. In particular, in genus 0, they satisfy the Witten--Dijkgraaf--Verlinde--Verlinde (WDVV) equation.
These classes are convenient for defining the quantum product on
either the Chow ring or cohomology ring.  The associativity of that
product is a consequence of the WDVV equation.

A fundamental ingredient of Gromov--Witten theory is a description
of the ``boundary" of the moduli stack of maps. The locus of nodes
$\Sigma$ on the universal curve of the stack of stable maps
$\ocM_{0,n}(X,\beta)$ is a partial normalization of this boundary,
and has the beautiful description
$$\Sigma\ \  =\ \  \coprod_{\substack{A\sqcup B = \{1,\ldots, n\}\\\beta_1+\beta_2\ \  =\ \  \beta}}  \ocM_{0,A \sqcup  \bigstar}(X,\beta_1)\ \mathop\times\limits_X \ \ocM_{0,B \sqcup \bullet}(X,\beta_2). $$ Here the fibered product is taken with respect to the \emph{evaluation maps} $e_{\bigstar}:  \ocM_{0,A \sqcup\bigstar}(X,\beta_1)\to X$ and $e_{\bullet}:  \ocM_{0,B \sqcup\bullet}(X,\beta)\to X$.

At the bottom of this is the fact that if $C = C_1 \sqcup_p C_2$ is
a nodal curve with node $p$ separating it into two subcurves $C_1$
and $C_2$, then $C$ is a coproduct of $C_1$ and $C_2$ over $p$. It
follows from the universal property of a coproduct that
$$\Hom(C, X) \ \ =\ \  \Hom(C_1,X) \mathop\times\limits_{\Hom(p,X)} \Hom(C_2,X).$$ Since $   \Hom(p,X) = X$ we get the familiar formula
$$\Hom(C, X) \ \ = \ \ \Hom(C_1,X) \ \mathop\times\limits_{X}\  \Hom(C_2,X).$$ The decomposition of the boundary is immediate from this.

\subsection{Boundary of moduli - orbifold case}
When analyzing the orbifold case something new happens. The source
curve is a twisted curve $\cC = \cC_1 \sqcup_\cG \cC_2$ with node
$\cG$ which is a gerbe banded by $\bmu_r$ for some $r$. We show in Proposition \ref{Prop:gluing} that
the coproduct $\cC_1\sqcup^\cG \cC_2$ of $\cC_1$ and $\cC_2$ over $\cG$ exists, and it is immediate that $\cC_1\sqcup^\cG \cC_2\to \cC$ is an isomorphism. It follows
again that
$$\Hom(\cC, \cX) \ \ =\ \  \Hom(\cC_1,\cX)
\mathop\times\limits_{\Hom(\cG,\cX)} \Hom(\cC_2,\cX),$$ see Appendix
\ref{Sec:gluing}. Now the data of
a gerbe with a map to $\cX$ is not a point of $\cX$, but of a
fascinating gadget $\riX$ we call \emph{the rigidified cyclotomic
  inertia stack}, see Section \ref{Sec:riX}. (A different notation $\bar {\mathcal{X}}_1$
was used in \cite{AGV}, but we were convinced that the
present notation, as used by
Cadman\cite{Cadman1,Cadman2}, is more
appropriate). We get the formula
$$\Hom(\cC, \cX) \ \ = \ \ \Hom(\cC_1,\cX) \mathop\times\limits_{\riX}
\Hom(\cC_2,\cX).$$ The fibered product uses a natural morphism
$\Hom(\cC_2, \cX) \to \riX$ corresponding to $\cG \to  \cX$, and a {\em
  twisted} map $\Hom (\cC_1,\cX) \to \riX$ corresponding to $\cG \to
\cX$ \emph{with the band inverted}. This is necessary since the glued
curve $\cC$ is \emph{balanced}.
The resulting map of moduli stacks
$$ \coprod_{\substack{ A\sqcup B = \{1,\ldots, n\}\\  \beta_1+\beta_2 =
    \beta}}  \cK_{0,A
  \sqcup\checkbullet}(\cX,\beta_1)\mathop\times\limits_{\riX} \cK_{0,B
  \sqcup \bullet}(\cX,\beta_2) \ \ \lrar \ \   \cK_{0,A \sqcup
  B}(\cX,\beta)$$ is crucial for Gromov--Witten theory of stacks. Here
the fibered product is taken using an evaluation map $$e_\bullet:
\cK_{0,B \sqcup \bullet}(\cX,\beta_2) \ \ \lrar \ \  \riX$$ and a
   \emph{twisted evaluation map}
$$\check e_\checkbullet: \cK_{0,A \sqcup \checkbullet}(\cX,\beta_1) \
   \ \lrar \ \  \riX$$ necessary to make the glued curves balanced,
   see Section \ref{Sec:evaluation}.

\subsection{Gromov--Witten classes - orbifold case}

We can now define the orbifold Gromov--Witten classes by integrating
along evaluation maps. We use the formalism of Chow rings, though
any cohomology theory satisfying some reasonable assumptions works.

Since evaluation maps land naturally  in $\riX$, \emph{the correct
cohomological theory to use is not that of $\cX$ but rather of
$\riX$}. The fact that something like  $A^*(\riX)_\QQ$ or
$\H^*(\riX,\QQ)$ is of interest was recognized by physicists (see
e.g. \cite{DHVW},  \cite{Zaslow}) where 
 ``twisted
sectors" and ``orbifold Euler characteristics" were considered. Some
mathematical reasons were discussed in \cite{AGV}. Also, in
Kontsevich's remarkable messages to Borisov the same phenomenon
occurs. But all these are reasoned by analogy, ``delicious
reciprocal reflections, furtive caresses, inexplicable quarrels"
\cite{Weil}. Alas, in orbifold Gromov--Witten theory all this becomes
mundane once one understands that it
all follows from the fact that
$\cC$ is a coproduct.

Now, given $\gamma_1,\ldots,\gamma_n \in A^*(\riX)_\QQ$ we can
consider the class $$\check e_{n+1\ *}(e_1^*\gamma_1 \cup\ldots \cup
e_n^*\gamma_n\cap [\cK_{g,n+1}(\cX,\beta)]^\vir)\in A^*(\riX)_\QQ,$$
where $e_i :\cK_{g,n+1}(\cX , \beta) \to \riX$ are the evaluation
maps of Section \ref{Sec:evaluation}.  It turns
out that in the orbifold theory we need to multiply
this by the function $r:\riX \to \ZZ$ describing the index of the
gerbe (but see the second part of Proposition
\ref{Prop:lifted-evaluation} for a way out
of that annoyance). We thus define the Gromov--Witten
classes to be
$$\langle \gamma_1,\ldots,\gamma_n ,*\rangle_{g,\beta} \ \ = \ \ r
\cdot \check e_{n+1\ *}(e_1^*\gamma_1 \cup\ldots \cup
e_n^*\gamma_n\cap [\cK_{g,n+1}(\cX,\beta)]^\vir).$$  With this
definition Gromov--Witten theory goes through almost as in the
manifold case, though the proofs require some interesting changes. The
main result is Theorem \ref{Thm:WDVV}, in which the WDVV equation is
proven.

\subsection{Acknowledgements}The authors are grateful to L.
Borisov, B. Fantechi, T. Kimura, R. Kaufmann, M. Kontsevich, M. Olsson, and H. Tseng for
helpful conversations. The two referees checked the paper very thoroughly, which is much appreciated by the authors.

\section{Chow rings, cohomology and homology of stacks}\label{classical}

\subsection{Intersection theory on Deligne--Mumford stacks}

Throughout the paper we work over a fixed base field $k$ of
characteristic $0$. Let $\cX$ be a separated Deligne--Mumford stack
of finite type over $k$, $\pi: \cX \to X$ its moduli space.

The rational Chow group $A_{*}(\cX)_{\QQ}$ is defined as in
\cite{Mumford}, \cite{Gillet} and \cite{Vistoli}. One defines the
group of cycles on $\cX$ as the free abelian group on closed
integral substacks of $\cX$, then divides by rational equivalence. There is also an integral version of the theory, developed in \cite{EG} for quotient stacks and in \cite{Kresch} in general, but we will not need it.

These groups are covariant for proper morphisms of Deligne--Mumford
stacks. The pushforward $\pi_{*}: A_{*}(\cX)_{\QQ} \to
A_{*}(X)_{\QQ}$ is given by the following formula. Let $\cV$ be a
closed integral substack of $\cX$, and call $V$ its moduli space;
this is an integral scheme of finite type over $k$. Because of the
hypothesis on the characteristic, the natural morphism $V \to X$ is
a closed embedding, and we have
   \[
   \pi_{*}[\cV] = \frac{1}{r}[V],
   \]
where $r$ is the order of the stabilizer of a generic geometric point of $\cV$.

The homomorphism $\pi_{*}$  is an isomorphism. In what follows we
will always identify $A_{*}(\cX)_{\QQ}$ and $A_{*}(X)_{\QQ}$ via
$\pi_{*}$.

If $\cX$ is proper, we denote by
   \[
   \int_{\cX}: A_{*}(\cX)_{\QQ}  \to \QQ
   \]
the pushforward $A_{*}(\cX)_{\QQ} \to A_{*}(\Spec k)_{\QQ} = \QQ$

If $f: \cX \to \cY$ is an {l.c.i.}~morphism of constant codimension
$k$, then for any morphism $\cY' \to \cY$ we have a \emph{Gysin
homomorphism} $f^{!}:A_{*}(\cY')_{\QQ} \to A_{*}(\cX\times_{\cY}\cY')_{\QQ}$, of degree $-k$. These commute with proper pushfowards, and among
themselves (\cite{Vistoli}).

As in \cite{Fulton}, Chapter~17, with every such stack $\cX$ we can associate a
bivariant ring $A^{*}(\cX)_{\QQ}$, that gives a contravariant
functor from the 2-category of algebraic stacks of
finite type over $k$ to the category of commutative rings. The
product of two classes $\alpha$ and $\beta$ in $A^{*}(\cX)_{\QQ}$
will be denoted by $\alpha\beta$, or by $\alpha\cup \beta$. By
definition, $A_{*}(\cX)_{\QQ}$ is a module over $A^{*}(\cX)_{\QQ}$;
we indicate the result of the action of a bivariant class $\alpha
\in A^{i}(\cX)_{\QQ}$ on a class of cycles $\xi\in A_{j}(\cX)_{\QQ}$
as a cap product $\alpha\cap \xi\in A_{j-i}(X)_{\QQ}$. The pullback
ring homomorphism $A^{*}(X)_{\QQ} \to A^{*}(\cX)_{\QQ}$ is also an
isomorphism. The projection formula holds: if $f : \cX \to \cY$ is a
proper homomorphism, $\beta \in A^{*}(\cY)_{\QQ}$ and $\xi \in
A_{*}(\cX)_{\QQ}$, then
   \[
   f_{*}(f^{*}\beta \cap \xi) = \beta\cap f_{*}\xi.
   \]

If $\cE$ is a vector bundle on $\cX$, then there are Chern classes
$c_{i}(\cE) \in A^{i}(\cX)$, satisfying the usual formal properties.

Suppose that $\cX$ is smooth: then the homomorphism
$A^{*}(\cX)_{\QQ} \to A_{*}(\cX)_{\QQ}$ defined by $\alpha \mapsto
\alpha \cap [\cX]$ is an isomorphism. In this case we can identify
$A^{*}(\cX)_{\QQ}$ with $A_{*}(\cX)_{\QQ}$.

The moduli space $X$ will not be smooth. However, assuming for the
moment that $\cX$ is connected, hence irreducible, and we call $r$
the order of the automorphism group of a generic geometric point of
$\cX$, then by the projection formula we have
   \[
   \pi_{*}(\pi^{*}\alpha \cap [\cX]) = \alpha\cap \pi_{*}[\cX] =
      \frac{1}{r}\alpha \cap [X]
   \]
for any $\alpha\in A^{*}(X)_{\QQ}$; hence the homomorphism
$A^{*}(X)_{\QQ} \to A_{*}(X)_{\QQ}$ defined by $\alpha \mapsto
\alpha \cap [X]$ is also an isomorphism. The same holds without
assuming that $\cX$ is connected, because if $\cX_{i}$ are the
connected components of $\cX$ and $X_{i}$ is the moduli space of
$\cX_{i}$, then $X = \coprod_{i} X_{i}$. Hence $A_{*}(X)_{\QQ}$
inherits a ring structure from that of $A^{*}(X)_{\QQ}$, even though
$X$ is in general singular.

However, one should be careful: the isomorphism $\pi_{*}:
A_{*}(\cX)_{\QQ} \to A_{*}(X)_{\QQ}$ is not a homomorphism of rings,
unless $\cX$ is generically a scheme, because in general the
identity $[\cX]$ of $A_{*}(\cX)$ is not carried into the identity
$[X]$ of $A_{*}(X)$.

Suppose that $\cX$ is smooth and proper: then $X$ is a complete
variety with quotient singularities. We say that an element of  the
group $A_{1}(X)_{\QQ}$ of 1-dimensional cycles is \emph{numerically
equivalent to $0$} if $\int_{X}\alpha\cup \xi = 0$ for all $\alpha
\in A^{1}(X)_{\QQ}$. The elements of $A_{1}(X)$ whose images in
$A_{1}(X)_{\QQ}$ are numerically equivalent to $0$ form a subgroup;
we denote by $N(X)$ the quotient group. This is finitely generated.
Furthermore we denote by $N^{+}(X)$ the submonoid of
$N(X)$ consisting of effective cycles.

Let $\cE$ be a vector bundle on $\cX$. If $\xi\in A_{1}(X)$ is an
integral 1-dimensional class, we denote
   \[
   c_{1}(\cE)\cdot \xi := \int_{\cX}c_1(\cE)\cap\xi',
   \]
where $\xi'$ is the class in $A_{1}(\cX)_{\QQ}$ such that
$\pi_{*}\xi'$ equals the image of $\xi$ in $A_{1}(X)_{\QQ}$.

Notice the following fact. If $\alpha$ is the class in
$A^{1}(X)_{\QQ}$ such that $\pi^{*}\alpha = c_{1}(\cE)$, then
   \begin{align*}
      c_{1}(\cE)\cdot \xi &= \int_{\cX}c_1(\cE)\cap\xi'\\
      &= \int_{X}\pi_{*}(c_1(\cE)\cap\xi')\\
      & = \int_{X}\alpha\cap\xi;
   \end{align*}
hence $c_{1}(\cE)\cdot \xi$ only depends on the class of $\xi$ in
$N(X)$. This allows us define the rational number $c_{1}(\cE)\cdot
\beta$ for a class $\beta\in N(X)$. It is easy to see  that the denominators
in  $c_{1}(\cE)\cdot \beta$ are uniformly bounded, but the following
proposition gives a natural bound:

\begin{proposition}\label{Prop:bound-denominators}
Assume that $\cX$ is proper. For each geometric point $p: \Spec
\overline{k} \to \cX$ denote by $e_{p}$ the exponent of the
automorphism group of $p$, and call $e$ the least common multiple of
the $e_{p}$ for all geometric points of $\cX$. Then
   \[
   c_{1}(\cE)\cdot \beta \in \frac{1}{e}\ZZ
   \]
for any vector bundle $\cE$ on $\cX$ and any $\beta\in N(X)$.
\end{proposition}

\begin{proof}
First note that $c_{1}(\cE)\cdot \beta \in \ZZ$  when $\cE = \pi^*\cM$ for a bundle $\cM$ on $X$. Indeed, in this case
\begin{align*}c_{1}(\cE)\cdot \beta &= \int_\cX \pi^*c_{1}(\cM)\cap \beta' \\
&=  \int_X c_{1}(\cM)\cap \beta\in \ZZ.
\end{align*}

We may substitute $\cE$ with its determinant, and assume that $\cE$
is a line bundle.

The following is a standard fact; we include a proof below for lack of a suitable reference.

\begin{lemma}
Let $\cL$ be a line bundle on $\cX$.  The line bundle $\cL^{\otimes
e}$ on $\cX$ is the pullback of a line bundle $\cM$ on $X$.
\end{lemma}

To conclude the proof of the proposition,
   \begin{align*}
      c_{1}(\cE)\cdot \beta &=  \frac{1}{e}c_{1}(\cL^{\otimes e})\cdot \beta\\
         &=  \frac{1}{e}\pi^* c_{1}(\cM)\cdot \beta \in  \frac{1}{e}\ZZ
   \end{align*}
as required. \qed

\begin{proof}[Proof of the lemma.] Observe that this is equivalent to the statement
that $\pi_{*}\cL^{\otimes e}$ is a line bundle on $X$, and the
adjunction homomorphism $\pi^{*}\pi_{*}\cL^{\otimes e} \to
\cL^{\otimes e}$ is an isomorphism. This is a local statement in
the \'etale topology.

Let $p:\Spec \overline{k} \to \cX$ be a geometric point, $G_{p}$ its
automorphism group. By definition, the exponent of
$G_{p}$ divides $e$. The action of $G_{p}$ on the fiber $\cL_{p}$ of
$\cL$ at $p$ is given by a 1-dimensional character $\chi: G_{p} \to
\overline{k}^{*}$, and therefore $\chi^{e}$ is the trivial character.
This implies that the action of $G_{p}$ on the fiber of
$\cL^{\otimes e}$ is trivial.

There is an \'etale neighborhood $\Spec\overline{k} \to U \to X$ of the image of $p$ in $X$, such that the pullback $U \times_{X}\cX$ is isomorphic to the quotient $[V/G_{p}]$, where $V$ is a scheme on which $G_{p}$ acts with a fixed geometric point $q: \Spec \overline{k} \to V$, with an invariant morphism $V \to U$ mapping $q$ to $p$; then the pullback $\cL_{[V/G_{p}]}$ corresponds to a $G_{p}$-equivariant locally free sheaf $\cL_{V}$ on $V$. We may assume that $U$ is affine. Then $V$ is also affine; since the characteristic of the base field is $0$, we can take an invariant  non-zero element the in fiber of $\cL^{\otimes e}_{V}$, and extend it to an invariant section of $\cL^{\otimes e}_{V}$. By restricting $V$ we may also assume that this section does not vanish anywhere: and then $\cL^{\otimes e}_{V}$ is trivial as a $G_{p}$-equivariant line bundle. This implies that the restriction $\cL^{\otimes e}_{[V/G_{p}]}$ is trivial, and so $\cL^{\otimes e}_{[V/G_{p}]}$ is the pullback of the trivial line bundle on $U$. This concludes the proof of the  Lemma and of the Proposition.
\end{proof}\noqed
\end{proof}

\subsection{Homology and cohomology of smooth Deligne--Mum\-ford stacks}

In this section the base field will be the field $\CC$ of complex
numbers. If $\cX$ is an algebraic stack of finite type over $\CC$,
then we can define the classical homology and cohomology of $\cX$ as
the homology and cohomology of the simplicial scheme associated with
a smooth presentation of $\cX$. More precisely, if $X_{1} \double X_{0}$ is a smooth presentation of $\cX$, we can obtain from it a simplicial scheme $X_{\bullet}$ in the usual fashion. To this we associate a simplicial space by taking the classical topology on each $X_{i}$. The homotopy type of the realization of this simplicial space is by definition the homotopy type of $\cX$, and its homology and cohomology are the homology and cohomology of $\cX$.

Alternatively, one can define the
classical site of $\cX$, and define cohomology as the
sheaf-theoretic cohomology of a constant sheaf on this site.
Homology can be defined by duality, as usual.

In this paper we are only going to need the cohomology and homology
with \emph{rational} coefficients of a \emph{proper} Deligne--Mumford
stack (mostly in the smooth case)
$\cX$. In this case a much more elementary approach is available.

Let $\cX$ be a separated Deligne--Mumford stack of finite type over
$\CC$, and let $\pi: \cX \to X$ be its moduli space. We define the homology
and cohomology $\H_{*}(\cX,\QQ)$ and $\H^{*}(\cX,\QQ)$ as
$\H_{*}(X,\QQ)$ and $\H^{*}(X,\QQ)$, respectively.

If $f: \cX \to \cY$ is a morphism of separated Deligne--Mumford
stacks of finite type over $\CC$ with moduli spaces $\pi: \cX \to X$
and $\rho: \cY \to Y$, this induces a morphism $g: X \to Y$ of
algebraic varieties over $\CC$. We define the pushforward $f_{*}:
\H_{*}(\cX,\QQ) \to \H_{*}(\cY,\QQ)$ and the pullback $f^{*}:
\H^{*}(\cX,\QQ) \to \H^{*}(\cY,\QQ)$ as
   \[
   g_{*}: \H_{*}(X,\QQ) \to \H_{*}(Y,\QQ)
   \quad\text{and}\quad
   g^{*}: \H^{*}(X,\QQ) \to \H^{*}(Y,\QQ)
   \]
respectively.

In particular, the pushfoward
   \[
   \pi_{*}: \H_{*}(\cX,\QQ) \to \H_{*}(X,\QQ),
   \]
is the identity function $\H_{*}(\cX,\QQ) \to \H_{*}(X,\QQ)$, and the pullback
   \[
   \pi^{*}: \H^{*}(X,\QQ) \to \H^{*}(\cX,\QQ)
   \]
is the identity $\H^{*}(X,\QQ) \to \H^{*}(\cX,\QQ)$.

With these definitions, $\H^{*}$ becomes a covariant functor from the
2-category of separated stacks of finite type over $\CC$ to the
category of graded abelian $\QQ$-vector spaces. Similarly, $\H^{*}$
becomes a contravariant functor from the same 2-category to the
category of graded-commutative $\QQ$-algebras.

Also, the cap product $\cap: \H^{*}(X, \QQ) \otimes \H_{*}(X,\QQ) \to
\H_{*}(X,\QQ)$ can be interpreted as a cap product $\cap: \H^{*}(\cX,
\QQ) \otimes \H_{*}(\cX,\QQ) \to \H_{*}(\cX,\QQ)$. If $f: \cX \to
\cY$ is a morphism of stacks, the projection formula
   \[
   f_{*}(f^{*}\alpha \cap \xi) = \alpha\cap f_{*}\xi
   \]
for any $\alpha\in \H^{*}(\cY,\QQ)$ and any $\xi\in \H_{*}(\cY, \QQ)$, holds.

Now assume that $\cX$ is proper. We define the cycle homomorphism
   \[
   \cyc_{\cX}: A_{*}(\cX)_{\QQ} \to \H_{*}(\cX,\QQ)
   \]
as the composition
   \[
   A_{*}(\cX)_{\QQ} \overset{\pi_{*}}{\longrightarrow}
   A_{*}(X)_{\QQ} \overset{\cyc_{X}}{\longrightarrow}
   \H_{*}(X,\QQ).
   \]
It is easy to see that the cycle homomorphism gives a natural
transformation of functors from the 2-category of proper
Deligne--Mumford stacks over $\CC$ to the category of graded
$\QQ$-vector spaces. If $\cV$ is a closed substack of $\cX$, this
 has a fundamental class $[\cV]$ in $A_{*}(\cX)_{\QQ}$; we
denote its image in $\H_{*}(\cX,\QQ)$ also by $[\cV]$, and call it
the \emph{homology fundamental class} of $\cV$.

Now assume that $\cX$ is smooth and proper. Then $X$ is a variety
with quotient singularities, hence it is a rational homology
manifold; thus we have Poincar\'e duality, that is, the homomorphism
   \[
   \PD_{X} := -\cap[X] : \H^{*}(X,\QQ) \longrightarrow \H_{*}(X, \QQ)
   \]
is an isomorphism. From this we get Poincar\'e duality on $\cX$;
however, one should be a little careful here, because the
fundamental class $[\cX]$, that we want to use to define Poincar\'e
duality on $\cX$, does not coincide with the fundamental class
$[X]$. In fact, if $\cX_{i}$ are the connected components of $\cX$,
then their moduli spaces $X_{i}$ are the connected components of $X$.
We have that $\pi_{*}[\cX_{i}] = \frac{1}{r_{i}}[X_{i}]$, if $r_{i}$
is the order of the automorphism of a generic geometric point of
$\cX_{i}$; hence we get the formula
   \[
   [\cX] = \sum_{i}\frac{1}{r_{i}}[X_{i}]
   \]
in $A_{*}(\cX)_{\QQ}$, and hence also in $H_{*}(\cX,\QQ)$.
We define the Poincar\'e duality homomorphism on $\cX$ as
   \[
   \PD_{\cX} := -\cap[\cX] : \H^{*}(\cX,\QQ) \longrightarrow \H_{*}(\cX, \QQ);
   \]
this is an isomorphism, because it is an isomorphism for every
connected component of $\cX$.

\section{The cyclotomic inertia stack and its rigidification} \label{Sec:riX}

\subsection{Cyclotomic inertia}
Let $\cX$ be a finite type Deligne--Mumford stack over $k$.
\begin{definition} We define a category $\iX[r]$, fibered over the
category of schemes, as follows:
\begin{enumerate}
\item Objects $\iX[r](T)$ consist of pairs $(\xi,\alpha)$ where $\xi$ is an object of $\cX$ over $T$,  and
$$\alpha: (\bmu_r)_T \to \cAut_T(\xi)$$ is an injective morphism of
group-schemes.  Here $(\bmu_r)_T$ is $\bmu_r \times T$.
\item An arrow from $(\xi,\alpha)$ over $T$ to $(\xi',\alpha')$ over $T'$ is an arrow $F: \xi \to \xi'$ making the following diagram commutative:
$$\xymatrix{
(\bmu_r)_T \ar[r]\ar_\alpha[d] & (\bmu_r)_{T'} \ar^{\alpha'}[d]\\
 \cAut_T(\xi) \ar[r]\ar[d]           & \cAut_{T'}(\xi') \ar[d]\\
                               T \ar[r]                  &               T'
}$$
where $(\bmu_r)_T  \to   (\bmu_r)_{T'} $ is the projection and $ \cAut_T(\xi) \to \cAut_{T'}(\xi') $ is the map induced by $F$.
\end{enumerate}
\end{definition}

It is evident that this category is fibered in groupoids over the category of schemes. There is an obvious morphism $\iX[r]\to \cX$ which sends $(\xi,\alpha)$ to $\xi$.

\begin{proposition}\label{Prop:inertia-finite}\label{Prop:X1-rep+finite}
The category $\iX[r]$ is a Deligne--Mumford stack, and the functor $\iX[r]\to \cX$ is representable and finite.
\end{proposition}

\begin{proof}
We verify that $\iX[r]\to \cX$ is representable and finite, which implies that $\iX[r]$ is a Deligne--Mumford stack. Consider a morphism $T \to \cX$ corresponding to an object $\xi$. The fibered product $\iX[r] \times_\cX T$ is an open and closed subscheme of the finite $T$-scheme
$\Hom_T( (\bmu_r)_T, \cAut_T(\xi)).$
\end{proof}

\begin{proposition}
Given an isomorphism  $\bmu_r \stackrel{\phi}{\lrar} \ZZ/r\ZZ$ there is an induced isomorphism $\iX[r]\simeq \cI(\cX,r)$, where $ \cI(\cX,r) \hookrightarrow  \cI(\cX)$ is the open and closed substack of the inertia stack of $\cX$, consisting of pairs $(\xi,g)$ with $g\in \cAut(\xi)$ of order $r$ at each point.
\end{proposition}

\begin{proof}
The data of an element $g\in \cAut(\xi)$ over $T$ of order $r$ at each point is equivalent to an injective group-scheme homomorphism $(\ZZ/r\ZZ)_T \to \cAut_T(\xi)$. Composing with $\bmu_r \stackrel{\phi}{\lrar} \ZZ/r\ZZ$ we get the result.
\end{proof}

\begin{corollary}\label{Cor:X1-smooth}
When $\cX$ is smooth, $\iX[r]$ is smooth as well.
\end{corollary}

\begin{proof}
It suffices to check this after extension of base field, so we may assume there exists an isomorphism $\bmu_r \stackrel{\phi}{\lrar} \ZZ/r\ZZ$, and by the Proposition it suffices to check the result for $\cI(\cX)$, which is well known.
\end{proof}

\begin{definition}
We define $\iX\  \ =\ \ \sqcup_r\ \iX[r],$ and we name it \emph{the cyclotomic inertia stack} of $\cX$.
\end{definition}

Note that, since $\cX$ is of finite type, $\iX[r]$ is empty except
for finitely many $r$, so $\iX$ is also of finite type, and in fact
finite over $\cX$ by Proposition \ref{Prop:inertia-finite}.

\subsection{Alternative description of cyclotomic inertia}
There is another less evident description of $\iX[r]$, which we give in the following definition.

\begin{definition}
We define a category $\iX[r]'$ over the category of schemes as follows.
\begin{enumerate}
\item Objects over a scheme $T$ consist of representable morphisms $\phi:(\cB\bmu_r)_T \to \cX$.
\item An arrow $\phi \to \phi'$ over $f: T \to T'$ is a 2-morphism $\rho: \phi\to  \phi'\circ f_*$ making the following diagram commutative:

$$\xymatrix{ (\cB\bmu_r)_T \ar^{f_*}[rr]\ar[rd]_\phi &&(\cB\bmu_r)_{T'} \ar[ld]^{\phi'} \\
& \cX
}$$
\end{enumerate}
\end{definition}

It is clear that this category is fibered in groupoids over the category of schemes.

\begin{definition}
We define a morphism  of fibered categories
$$\iX[r]' \to \iX[r]$$ as follows:
\begin{enumerate}
\item Given an object $$\xymatrix{
\cB(\bmu_r)_T \ar^\phi[r]\ar[d] & \cX \\ T
}$$ we obtain a pair $(\xi, \alpha)$ as follows:  $\xi$ is obtained by composing $\phi$ with the section $T\to \cB(\bmu_r)_T$ associated with the trivial $\bmu_r$-torsor $\bone_{\bmu_r,T}$. The homomorphism $\alpha$ is  the associated map of automorphisms
 $$ (\bmu_r)_T = \cAut_T(\bone_{\bmu_r,T}) \lrar \cAut_T(\xi),$$ which is injective since $\phi$ is representable.
 \item Given an arrow $\rho$ as above, we obtain $F: \phi(\bone_{\bmu_r,T})  \to \phi'(\bone_{\bmu_r,T'})$ by completing the following diagram:

 $$\xymatrix{
\phi(\bone_{\bmu_r,T}) \ar[d]_\rho\ar@{..>}[rrdd]^F \\
(\phi'\circ f_*)(\bone_{\bmu_r,T})\ar@{=}[d]\\
\phi'\left((\bone_{\bmu_r,T'})_T\right)\ar[rr] &&\phi'(\bone_{\bmu_r,T'})
}$$

\end{enumerate}
\end{definition}

\begin{proposition}
The morphism $\iX[r]' \to \iX[r]$ is an equivalence of fibered categories.
\end{proposition}

\begin{proof}
It is enough to show that for each scheme $T$, the induced functor on the fiber $\iX'(r)(T) \to \iX[r](T)$ is an equivalence. In the following proof we will write the action of a group on a torsor on the left, not on the right, as is more customary.

 {\sc Step 1: the functor is faithful.} Given two objects $\phi, \phi':\cB(\bmu_r)_T \to \cX$, suppose that $\alpha: \phi \to \phi'$ is a 2-arrow. We need to show that for any $T$-scheme $f: U \to T$ and any $\bmu_{r}$-torsor $P \to U$, the arrow $\alpha_{P\to U}: \phi(P\to U) \to \phi'(P\to U)$ is uniquely determined by $\alpha_{\bone_{\bmu_{r},T}}: \phi(\bone_{\bmu_{r},T}) \to \phi'(\bone_{\bmu_{r},T})$. Let $\{U_{i} \to U\}$ be an \'etale covering, such that the pullbacks $P_{i} \to U_{i}$ are all trivial. Then the pullbacks of $\phi(P\to U)$ and $\phi'(P\to U)$ to $\cX(U_{i})$ are $\phi(P_{i}\to U_{i})$ and $\phi'(P_{i}\to U_{i})$, respectively, and, since $\cX$ is a stack, the restrictions of a morphism $\phi(P\to U) \to \phi'(P\to U)$ is determined by its restriction to $\phi(P_{i}\to U_{i}) \to \phi'(P_{i}\to U_{i})$. So we may assume that $P\to U$ is trivial. But then there is a cartesian arrow $P \to \bone_{\bmu_{r},T}$, and an induced diagram
   \[
   \xymatrix@C+15pt{
   \phi(P\to U) \ar[r]^{\alpha_{P\to U}}\ar[d]
      & \phi'(P\to U) \ar[d]\\
   \phi(\bone_{\bmu_{r},T})\ar[r]^{\alpha_{\bone_{\bmu_{r},T}}}
      & \phi'(\bone_{\bmu_{r},T})
   }
   \]
that proves what we want.

{\sc Step 2: the functor is fully faithful.} Assume that $\beta:
\phi(\bone_{\bmu_{r},T}) \to \phi'(\bone_{\bmu_{r},T})$ is an arrow
in $\cX(T)$, commuting with the actions of $\bmu_{r}$. First
consider the case of a trivial  $\bmu_{r}$-torsor $P\to U$  over a
$T$-scheme. Choose a trivialization that induces a cartesian arrow
$P \to \bone_{\bmu_{r},T}$. Then, by definition of cartesian arrow,
there is a unique dotted arrow in $\cX(U_{i})$ that we can insert in
the diagram
   \[
   \xymatrix@C+15pt{
   \phi(P\to U) \ar@{.>}[r]^{\alpha_{P\to U}}\ar[d]
      & \phi'(P\to U) \ar[d]\\
   \phi(\bone_{\bmu_{r},T})\ar[r]^{\beta}
      & \phi'(\bone_{\bmu_{r},T})
   }
   \]
making it commutative. This arrow $\alpha_{P\to U}$ is independent of the chosen
trivialization, because $\beta$ commutes with the actions of $\mu_{r}$, and two
trivializations differ by a morphism $U \to \bmu_{r}$.

If $P \to U$ is not necessarily trivial, choose a covering $\{U_{i}
\to U\}$ such that the pullbacks $P_{i} \to U_{i}$ are trivial. We
have arrows $\alpha_{P_{i} \to U_{i}}: \phi(P_{i} \to U_{i})\to
\phi'(P_{i} \to U_{i})$ in $\cX(U_{i})$, and their pullbacks to
$U_{i}\times_{U}U_{j}$ coincide; hence they glue together to given
an arrow $\alpha_{P \to U}: \phi(P \to U)\to \phi'(P \to U)$. It is
easy to see that $\alpha_{P \to U}$ does not depend on the covering,
and defines a 2-arrow $\phi\to \phi'$ whose image in $\iX'(r)$
coincides with $\beta$.

{\sc Step 3: the functor is essentially surjective.} Let there be given an object $(\xi,\alpha)$ of $\iX[r](T)$, and let us construct a morphism $\phi: \cB(\bmu_r)_T \to \cX$ of fibered categories, whose image in $\iX'(r)(T)$ is isomorphic to $(\xi,\alpha)$.

Let $P \to U$ be a $\bmu_{r}$-torsor, where $U$ is a $T$-scheme: we will define an object $\eta$ of $\cX(U)$, that is a twisted version of the pullback $\xi_{U}$ to $U$, by descent theory, using the action of $\bmu_{r}$ on $\xi$. The facts that we are going to use are all in \cite{descent}, Sections 3.8 and 4.4.

Consider the pullback $\xi_{P}$ of $\xi$ to $P$. The morphism $\bmu_{r}\times P \to P\times_{U}P$ defined as a natural transformation via Yoneda's Lemma by the usual rule $(\zeta, p) \mapsto (\zeta p,p)$ is an isomorphism. The pullbacks of $\xi$ to $P\times_{U}P = \bmu_{r}\times P$ along the first and second projection coincide with the pullback $\xi_{\bmu_{r}\times P}$. On the other end, the group scheme $\bmu_{r}$ acts on $\xi$, so the projection $\bmu_{r}\times P \to \bmu_{r}$ induces an automorphism of $\xi_{\bmu_{r}\times P}$, that gives descent data for $\xi_{P}$ along the \'etale covering $P \to U$. These descent data are effective, and define an object $\eta$ of $\cX(U)$. So we have assigned to every object of $\cB(\bmu_{r})_{T}(U)$ an object of $\cX(U)$; this is easily seen to extend to a morphism of fibered categories $\phi: \cB(\bmu_{r})_{T} \to \cX$.

Let  $P = (\bmu_{r})_{T} \to T$ be the trivial torsor. We claim that the object $\eta \eqdef \phi(P \to T)$ of $\cX(T)$ is isomorphic to $\xi$. In fact, the object with descent data defining $\xi$ is $\xi_{P}$, with the descent data given by the identity on $\xi_{P\times_{T}P}$. Then the projection $P \to \bmu_{r}$ defines an automorphism of $\xi_{P}$ in $\cX(P)$, that is easily seen to descend to an isomorphism $\xi \simeq \eta$ in $\cX(T)$. This isomorphism is $\bmu_{r}$-equivariant, because $\bmu_{r}$ is commutative, hence the image of $\phi$ in $\iX(T)$ is isomorphic to $(\xi,\alpha)$, as we wanted.
\end{proof}

\subsection{The stack of gerbes in $\cX$}

We introduce a stack $\riX$ closely related to $\iX$, which will
play an important role below.  It will be defined in terms of
morphisms of gerbes.  Recall that a gerbe over
a scheme $X$ is an fppf stack $F$ over
$X$ such that
\begin{itemize}
\item there exists an fppf covering $\{X_i \to X\}$ such that
$F(X_i)$ is not empty for any $i$, and

\item given two objects $a$ and $b$ of $F(T)$, where $T$ is an
$X$-scheme, there exists a covering $\{T_i \to T\}$ such that the
pullbacks $a_{T_i}$ and $b_{T_i}$ are isomorphic in $F(T_i)$.

\end{itemize}

Slightly more generally, a stack $F$ over a {\em stack} $\cX$  is a gerbe if for any morphism  $V \to \cX$ with $V$ a scheme, $F_V \to V$ is a gerbe.

If $G$ is a sheaf of abelian groups over $X$, we say that $F$ is
banded by $G$ if for each object $a$ of $F(T)$ we have an
isomorphism of sheaves of groups $Aut_T(a)$ with $G_T$. This should
be functorial, in the obvious sense.

If $G$ is not abelian, one can still define a gerbe banded by $G$,
but the definition is more subtle. If $G$ is a group scheme of
finite type, it is not hard to see that every gerbe banded by $G$ is
an algebraic stack.

We first need the following definition.

\begin{definition}
Define a 2-category $\tworiX[r]$ with a functor to the category of schemes as
follows:
\begin{enumerate}

\item An object over a scheme $T$ is a pair $(\cG,\phi)$, where  $\cG \to T$ is a gerbe banded by $\bmu_{r}$, and $\cG \stackrel{\phi}{\to} \cX$ is a representable morphism.
\item A morphism  $(F,\rho): (\cG,\phi) \to (\cG',\phi')$ consists of a morphism $F : \cG \to \cG'$ over some $f: T \to T'$, compatible with the bands,
 and a 2-morphism
$\rho: \phi\to  \phi'\circ F$ making the following diagram commutative:
$$\xymatrix{ \cG \ar^{F}[rr]\ar[rd]_\phi && \cG' \ar[ld]^{\phi'} \\
& \cX
}$$
\item A 2-arrow $(F,\rho) \to  (F_{1},\rho_{1})$ is a usual  2-arrow $\sigma: F \to F_{1}$ compatible with $\rho$ and $\rho_{1}$ in the sense that the following diagram is commutative:
   \[
   \xymatrix{ & \phi \ar[dl]_{{\rho}} \ar[dr]^{\rho_{1}}\\
    \phi'\circ F \ar[rr]^{\phi'(\sigma)} && \phi'\circ F_{1}
   }
   \]
\end{enumerate}
\end{definition}
\begin{lemma}
The 2-category   $\tworiX[r]$  is equivalent to a category.
\end{lemma}

\begin{proof}
Since all 2-arrows are isomorphisms, it suffices to show that the automorphism group of any  1-arrow is trivial. This is the content of the following general lemma.
\end{proof}

\begin{lemma} \label{lem:representable->1-category}
Suppose given a diagram
$$\xymatrix{ \cG \ar^{F}[rr]\ar[rd]_\phi && \cG' \ar[ld]^{\phi'} \\
& \cX }$$ where $\cG$, $\cG'$ and $\cX$ are categories  fibered in
groupoids over a base category, with a 2-arrow $\rho: \phi \to
\phi'\circ F$ making the diagram commutative.  Assume $\phi'$ is
faithful. Then an automorphism of $F$ compatible with $\rho$ is the
identity.
\end{lemma}

\begin{proof}
Take an object $\xi$ of $\cG$ over some $T$ in the base category. We get a diagram
  \[
   \xymatrix{ & \phi(\xi) \ar[dl]_{{\rho_{\xi}}} \ar[dr]^{\rho_{\xi}}\\
    \phi'( F(\xi)) \ar[rr]^{\phi'(\sigma_{\xi})} && \phi'( F(\xi))
   }
   \]
But $\phi'(\sigma_{\xi})$ lies over the identity $T \to T$. Since $\cX$ is fibered in groupoids, it follows that $ \phi'(\sigma_{\xi})$ is the identity. Since $\phi'$ is faithful, $\sigma_{\xi}$ is the identity, which is what we wanted.
\end{proof}

\begin{definition}\label{definition-riX}
We define $\riX[r]$ to be the category associated with the 2-category $\tworiX[r]$,
where arrows in  $\riX[r]$ are 2-isomorphism classes of 1-arrows in
$\tworiX[r]$.
\end{definition}
We note that $\riX[r]$ is a category fibered in groupoids over the
category of schemes.

\begin{remark}\label{rigidification-remark}
There is a tautological morphism $\iX[r] \to \riX[r]$, defined as follows. An object of $\iX[r](T)$ corresponds to a representable morphism $(\cB\bmu_{r})_{T} \to \cX$. Since $(\cB\bmu_{r})_{T}$ is a gerbe over $T$ banded by $\bmu_{r}$, this gives an object of $\riX$. An arrow in $\iX[r]$ gives an arrow in $\riX[r]$ in the obvious way.

We can also define a fibered category, an object of
which is an object of $\riX[r]$, together with a section of the
gerbe. There is an obvious forgetful functor into $\riX[r]$, which
exhibits it as the universal gerbe over $\riX$. Since a gerbe over
$T$ banded by $\bmu_{r}$ with a section has a canonical isomorphism
with $(\cB\bmu_{r})_{T}$, this universal gerbe is evidently
isomorphic to $\iX[r] \to  \riX$.
\end{remark}

In the next section we give an alternative description of $\riX[r]$,
which in particular  shows that it is a Deligne--Mumford stack.

\begin{definition}
We define $\riX\  \ =\ \ \sqcup_r\ \riX[r],$ and we name it \emph{the
stack of cyclotomic gerbes} in $\cX$.
\end{definition}

\subsection{The rigidified cyclotomic inertia stack}
Consider the stack $\iX[r]$. By definition, for every object
$(\xi,\alpha)$ of $\iX[r]\bigl(T\bigr)$ we have a canonical {\em
central} embedding $(\bmu_{r})_{T}
\stackrel{\iota_{(\xi,\alpha)}}{\hookrightarrow} \cAut_{T}
(\xi,\alpha)$. By \cite{ACV}, Theorem 5.1.5 (see also Appendix~\ref{Sec:rigidification}), there
exists a {\rm rigidification} denoted there $\iX[r] \to
\iX[r]^{\bmu_{r}}$. We are adopting the better notation proposed by
Romagny in \cite{Romagny}, and denote this by
$\iX[r]\thickslash\bmu_{r}$.

\begin{proposition}
We have an equivalence of fibered categories
$$\iX[r]\thickslash\bmu_{r} \to \riX[r]$$ so that the following
diagram is commutative:
   \[
   \xymatrix{
    & \iX[r]\ar[dl]\ar[dr]\\ \iX[r]\thickslash\bmu_{r}\ar[rr]&&\riX[r]
   }
   \]
The right diagonal arrow is described in Remark~\ref{rigidification-remark}.
\end{proposition}

\begin{proof}
We will use the modular interpretation
$\iX[r]\slashprime \bmu_{r}$ of the rigidified stack
$\iX[r]\thickslash \bmu_{r}$ given in
Proposition~\ref{prop:equivalence-rigidification}.  There is an
obvious morphism of fibered categories $\iX[r]\slashprime \bmu_{r}
\to \riX[r]$, defined as follows. An object of $\iX[r]\slashprime
\bmu_{r}(T)$ consists of a gerbe $\cG \to T$ banded by $\bmu_{r}$,
and a $\bmu_{r}$-2-equivariant morphism $\cG \to \iX[r]$, in the sense of Section~\ref{subsection-moduli-interpretation}: this can
be composed with the projection $\iX[r] \to \cX$ to get a
representable morphism $\cG \to \cX$. This function on objects can
be extended to a function on arrows, and it defines the desired
functor.

We claim that this an equivalence: let us construct an inverse $\riX
\to \iX\slashprime \bmu_{r}$. Consider an object of $\riX$,
consisting of a gerbe $\cG \to T$ banded by $\bmu_{r}$, and a
representable morphism $\phi: \cG \to \cX$. Given an object $\xi$ of
$\cG(U)$, where $U$ is a $T$-scheme, the morphism $\phi$ induces a
homomorphism of group-schemes $\alpha: H_{U} = \cAut_{U,\cG}(\xi)\to \cAut_{U,\cX}(\phi\xi)$; and this homomorphism is an
embedding, because a representable morphism is also faithful, so
$(\xi,\alpha)$ is an object of $\iX$. This function from objects of
$\cG$ to objects of $\iX$ extends naturally to a morphism $\cG \to
\iX$; and this morphism is $\bmu_{r}$-2-equivariant, by definition.

 We leave it to the reader to extend this
to a functor $\riX \to \iX\slashprime \bmu_{r}$, and show that it
gives an inverse to the functor $\iX[r]\slashprime \bmu_{r} \to
\riX[r]$ above. The commutativity of the diagram is straighforward.
\end{proof}

\begin{corollary}
If $\cX$ is smooth, $\riX$ is smooth as well. If $\cX$ is proper, $\riX$ is proper as well.
\end{corollary}

\begin{proof}
This follows from Proposition \ref{Prop:X1-rep+finite} and Corollary \ref{Cor:X1-smooth}, since the morphism $\iX \to \riX$ is proper, \'etale and surjective.
\end{proof}

\subsection{Changing the band by a group automorphism}\label{Sec:changing-band}
There is an involution $\iota:\riX\to \riX$ defined as follows:
given a gerbe $\cG \to T$ banded by a group-scheme  $G$ and an
automorphism $\tau: G \to G$, we can change the banding   of the
gerbe through the automorphism $\tau$. Applying this procedure to
the gerbe $\cG$ of an
 object $\cG \to \cX$ of $\riX[r](T)$, we get another object
 $\sideset{^\tau}{}{\mathop{\cG}}
 \to \cX$
of  $\riX[r](T)$. When $\tau:\bmu_{r}\to \bmu_{r}$ is the inversion
automorphism $\zeta \mapsto \zeta^{-1}$, this induces an involution
of $\riX[r]$. Applying this to each piece of $\riX$ separately, we
obtain the desired involution of $\iota:\riX\to \riX$.

\subsection{The tangent bundle lemma}
\begin{lemma}\label{tangent}  Let $S$ be a scheme, and let $f:S \to \riX$ be
a morphism.  Let
$$\xymatrix{\cG \ar[r]^F \ar[d]_\pi & \cX\\
S}$$
be the associated diagram.  Then there is a canonical isomorphism between
$\pi_*(F^*(T_\cX))$ and $f^*(T_{\riX})$.
\end{lemma}

\begin{proof}
Consider the universal gerbe $\iX \to \riX$ and the diagram of smooth stacks
  $$\xymatrix{\iX \ar^{F_1}[r]\ar^\varpi[d] & \cX \\ \riX . }
  $$
  Given a morphism $f:S \to \riX$ as in the lemma, we have a fiber diagram
   $$\xymatrix{\cG \ar^g[r]\ar@/^1.5pc/[rr]^{F} \ar_\pi[d]    &\iX \ar^{F_1}[r]\ar^\varpi[d] & \cX \\ S \ar^f[r] & \riX . }
  $$
  Since $\varpi$ is flat and $\varpi_*$ is exact on coherent sheaves,  for any locally free sheaf $\cH$ on $\iX$ we have
  $f^* \varpi_* \cH = \pi_*g^*\cH$. Therefore it suffices  to check that
  $\varpi_*(F_1^*(T_\cX))$ and $T_{\riX}$ are canonically isomorphic.

We have a natural morphism of sheaves $T_{\iX} \to F_1^* T_\cX$, giving a morphism
$\varpi_*T_{\iX} \to \varpi_* F_1^* T_\cX$. Since $\iX \to \riX$ is an \'etale gerbe, we have $\varpi_*T_{\iX} \cong T_{\riX}$, giving a morphism
$T_{\riX} \to \varpi_* F_1^* T_\cX$. We can check that this is an isomorphism by pulling back to geometric points.

Over a geometric point $y$ of $\riX$, we can identify the
fiber of $\varpi$ with $\cB\bmu_r$. This gives a lift of $y$ to $\iX$, and
so the point $y$ maps to a geometric
point $x$ of $\cX$, with stabilizer $G$.  We can locally describe
$\cX$ around $x$ as $[U/G]$. The pullback $T$ of the tangent space of $\cX$
to $y$ has a natural action of $\bmu_r$, and the fiber of
$\varpi_*(F_1^*T_\cX)$ at $y$ is naturally the space of invariants
$T^{\bmu_r}$.  Given a local chart of $\cX$ of the form $[U/G]$, the stack $\iX$ has a local chart
given by $[U^{\bmu_r}/C(\bmu_r)]$.  Since $T_{U^{\bmu_r}}=T_U^{\bmu_r}$, we obtain a natural isomorphism $T_{\riX,y} \cong T^{\bmu_r}$, which is what we need.
\end{proof}

\section{Twisted curves and their maps} \label{Sec:twisted-curves}

The foundation of the theory of stable maps to an orbifold rests on the notion of a \emph{twisted curve}. Over an algebraically closed field, a twisted curve is a connected, one-dimensional Deligne--Mumford stack which is \'etale locally a nodal curve, and which is a scheme outside the marked points and the singular locus.  Moreover, we will always include the condition that the nodes be balanced, that is, formally locally near a node, the stack is isomorphic to $$\big[\Spec (k[x,y]/(xy))\, \big/\, \bmu_r\big]$$ where the action of $\mu_r$ is given by $\zeta (x,y) = (\zeta \cdot x, \zeta^{-1} \cdot y)$.
In particular, the coarse moduli space of a twisted curve is always a nodal curve.

The notion of a family of twisted curves is straightforward, but involves one novelty.  Although the relative coarse moduli scheme of a family of twisted pointed curves is a family of prestable curves -- and hence comes with sections corresponding to the marked points -- a family of twisted curves need not have sections.  The marked point on each fiber instead gives rise to a \emph{gerbe} banded by $\bmu_r$ where $r$ is the order of the inertia group at the twisted point.

\subsection{The stack of twisted curves}\label{Sec:stack-twisted-curves}

If we define a groupoid $\fM^{\tw}_{g,n}$ whose $S$ points are given
by families of twisted curves over $S$, that is morphisms $\pi: \cC
\to S$ which are flat, together with a collection of $n$ disjoint
gerbes over $S$ embedded in $\cC$, such that the geometric fibers of
$\pi$ are $n$ pointed twisted curves of  genus $g$, then the
following results are proven in \cite{Olsson-twisted} .

\begin{enumerate}
\item $\fM^{\tw}_{g,n}$ is a smooth algebraic stack, locally of finite type.
\item If we bound the topological type, including the twisting at marked points and nodes, we get a stack of finite type.
\item A formal deformation space $\Delta^{\tw}$ of a twisted curve
  $\cC$ with nodes $q_{1},\ldots, q_{s}$ of indices
  $r_{1},\ldots,r_{s}$ can be obtained from a given formal deformation
  space $\Delta = \Spf R$ of the coarse curve $C$ as follows: Let
  $D_{i}$ be the divisor in $\Delta$ corresponding to $q_{i}$, with
  defining equations $f_{i}$. Then $$\Delta^{\tw} = \Spf R\llbracket x_{1},\ldots,x_{s}\rrbracket/\left(x_{1}^{r_{1}} - f_{1},\ldots x_{s}^{r_{s}} - f_{s}\right)$$ is a formal deformation space of $\cC$.
\end{enumerate}

\subsection{The smooth locus of a twisted pointed curve}
The theory in this section is due to the authors and to C.~Cadman
independently (\cite{Cadman1}, Section~2). The reader is advised to read Appendix \ref{roots} for the notion of root stacks.

Suppose that $C \to S$ is an $n$-pointed nodal curve; call $s_i
\colon S \to C$ the sections, $S_i \subseteq C$ their images,
$\sigma_i$ the canonical section of $\cO(S_i)$. Given positive
integers $d_1$, \dots,~$d_n$, we define a stack over $S$ as the
fibered product
   \[
    \cC[d_1, \dots, d_n]\ \  =\ \  \radice{d_1}{(\cO(S_1),
   \sigma_1)/C}\ \mathop\times\limits_C\ \cdots\ \mathop\times\limits_C\
   \radice{d_n}{(\cO(S_n), \sigma_n)/C}.
   \]
   See Appendix \ref{roots} for an explanation of the notation.
The locus where the projection
$$\radice{d_i}{(\cO(S_i), \sigma_i)/C} \to C$$ is not an isomorphism
coincides with $S_i$ when $d_{i} > 1$; if $d_{i} = 1$ the projection
is an isomorphism. Outside of the locus $S_i$
the stack $\cC[d_1, \dots, d_n]$ is isomorphic to $C$; over $S_i$ we
have an embedding $$\radice{d_i}{\cO_C(S_i)}\mid_{S_i} \into
\radice{d_i} {(\cO(S_i),\sigma_i)/C};$$ since we have that for $j
\neq i$ the morphism $\radice{d_j}{(\cO(S_j), \sigma_j)/C} \to C$ is
an isomorphism in a neighborhood of $S_i$, we also have an embedding
$$
\radice{d_i} {\cO_C(S_i)}\mid_{S_i} \into \cC[d_1, \dots, d_n].
$$
It
is easily checked that, after taking the
$\radice{d_i}{\cO_C(S_i)}\mid_{S_i}$ as markings, the stack
$\cC[d_1, \dots, d_n]$ is a twisted curve over $S$ with moduli space
equal to $C$.

The following theorem shows that the smooth part of a twisted curve
is uniquely characterized by its moduli space and its indices along
the markings. If $\cC \to S$ is a twisted curve, we will denote its
smooth part by $\cC_\mathrm{sm}$.

\begin{theorem}\label{thm:description-twistedcurves}
Let $\cC \to S$ be an $n$-pointed twisted curve with moduli space
$\pi: \cC \to C$. Assume that its index at the $i\th$-section is
constant for all $i$, and call it $d_i$. Denote by $\Sigma_i$ the
$i\th$ marking, by $s_i \colon S \to C$ the section corresponding to
$\Sigma_i$, and by $N_i$ the normal bundle to $S$ along $s_i$.

\begin{enumerate}

\item There is a canonical isomorphism of twisted curves
   \[
   \cC_\mathrm{sm} \simeq \cC[d_1, \dots, d_n]_\mathrm{sm}
   \]
inducing the identity on $C_\mathrm{sm}$.

\item There is a canonical isomorphism of gerbes over $S$
   \[
   \Sigma_i\simeq \radice{d_i}{N_i/S}.
   \]
In particular, $\Sigma_{i}$ is canonically   banded by $\bmu_{d_{i}}$.

\end{enumerate}

\end{theorem}

\begin{proof}
For part~(1), notice that we have an equality of a Cartier divisor $\pi^{*}S_{i}  = d_{i}\Sigma_{i}$ on $\cC$; this induces a morphism $\cC \to \cC[d_1, \dots, d_n]$. This is an isomorphism outside of the marked points and the nodes, and is easily seen to be representable outside of the nodes. To check that the restriction $\cC_\mathrm{sm} \to \cC[d_1, \dots, d_n]_\mathrm{sm}$ is an isomorphism it is enough to restrict to the geometric fibers of $S$, because $\cC$ and $\cC[d_1, \dots, d_n]$ are flat over $S$; but the morphism $\cC_\mathrm{sm} \to \cC[d_1, \dots, d_n]_\mathrm{sm}$ restricted to a geometric fiber is representable, finite and birational, and the stacks appearing are smooth, hence it is an isomorphism.

Part~(2) follows immediately from part~(1).
\end{proof}

\subsection{Twisted stable maps} 
Let $\cX$ be a Deligne--Mumford stack, $g,n$ non-negative integers
and $\beta$ a curve class on the coarse moduli space $X$. Associated
with
this data we have the stack $\cK_{g,n}(\cX,\beta)$ of $n$-pointed
twisted stable maps into $\cX$ of genus $g$ and class $\beta$.  This
classifies stable maps from $n$-pointed twisted genus $g$ curves to
$\cX$ of degree $\beta$.  Precisely, an object of this stack over a
scheme $T$ consists of the data of a family of twisted curves $\cC
\to T$, $n$ gerbes $\Sigma_i \subset \cC$, and a representable
morphism $\cC \to \cX$, such that the induced maps of underlying
coarse moduli spaces give a family of $n$ pointed genus $g$ stable
maps to $X$.  We refer the reader to \cite{AV} for a construction
and a more detailed discussion of this stack.

Let $\cC \to T$ be a twisted $n$-pointed curve with $T$ connected,
and let $1\leq i \leq n$. Assume that the index of the $i$-th
marking is the integer $r$. Note that, by Theorem
\ref{thm:description-twistedcurves}, part (2) we have that
$\Sigma_{i}^{\cC}$ is canonically banded by $\bmu_{r}$.

\subsection{Evaluation maps}\label{Sec:evaluation}\hfil
\begin{definition}\hfil
\begin{enumerate}
\item Assume given a twisted stable map $f: \cC \to \cX$ over a base $T$. We define
$$e_{i}(f)\in  \riX[r]\bigl(T\bigr)$$
to be the object  associated with the diagram
   \[
   \xymatrix{
   \Sigma_{i}^{\cC} \ar[r]^{f\mid_{ \Sigma_{i}^{\cC}}} \ar[d] & \cX \\
   T
   }
   \]
By Definition~\ref{definition-riX}, this defines a morphism
$$e_{i}: \cK_{g,n}(\cX,\beta) \to \riX,$$
which we call \emph{the $i$-th evaluation map}. \item The morphism
$\check e_{i} := \iota\circ e_{i}$, where $\iota:\riX\to \riX$ is
the involution defined in section \ref{Sec:changing-band}, is called
\emph{the $i$-th twisted evaluation map}.
\end{enumerate}
\end{definition}

\subsection{The virtual fundamental class}
The key technical point in developing Gromov-Witten theory for $\cX$ is the
construction of the \emph{virtual fundamental class} $[\cK_{g,n}(\cX,\beta)]^\vir$ in $A_*(\cK_{g,n}(\cX,\beta))$.  By \cite{BF}
and \cite{LT}, what is needed to construct this class is a perfect obstruction theory on this moduli stack.  Following the methods of \cite{BF}, we will mean by this a morphism in the derived category
$$\phi:E \to \bL_{\cK_{g,n}(\cX,\beta)/\fM^\tw_{g,n}}$$
such that
\begin{itemize}
\item $E$ is locally equivalent to a two term complex of locally free sheaves, and
\item $\H^0(\phi)$ is an isomorphism and $\H^{-1}(\phi)$ is surjective.
\end{itemize}
As in the case of ordinary stable maps, there is a natural perfect
obstruction theory with $E = \bR\pi_*(f^*T\cX)^\vee$.  The proof of
this is identical to the proof for
ordinary stable maps since what is needed are formal properties of
the cotangent complex and Illusie's results \cite{Illusie} relating
these to the deformation theory of morphisms.  Since the theory of
the cotangent complex for Artin stacks has been developed in \cite{LMB} and corrected
in \cite{Olsson-sheaves}, and since
Illusie explicitly works in the general setting of ringed topoi, all
the necessary generalizations have already been established.

Specifically, in the
discussion of \cite{Behrend} page 604, immediately after Proposition
4, one relies on the claim that $$\phi:\bR\pi_*(f^*T\cX)^\vee \to
\bL_{\cK_{g,n}(\cX,\beta)/\fM^\tw_{g,n}}$$ is a perfect {\em
relative} obstruction theory. This relative case, discussed in
section 7 of \cite{BF} (page 84 onward), reworks the absolute case
discussed earlier in that paper. The crucial result in \cite{BF} is
Proposition 6.3, where $\cC$ is assumed to be a Gorenstein and {\em
projective} curve. As explained above, projectivity is not necessary
for deformation theory (i.e.  \cite{BF} Theorem 4.5) - it works just
as well for a proper Deligne--Mumford stack.  Both in  \cite{BF}
Proposition 6.3 and in \cite{BF} Lemma 6.1 on which it relies, one
also needs relative duality, which is ``well known" for proper
Gorenstein Deligne--Mumford stacks; for a twisted curve $\cC$ with a
projective coarse moduli space it can be shown using a finite flat
Galois covering $D \to \cC$, ramified over auxiliary sections, which
can be constructed locally over the base.

We remark that one additional
feature of $E$ that is required in \cite{BF} is that $E$ admit a
global resolution (see discussion before \cite{BF}, Proposition 5.2).
Kresch's work on intersection theory for Artin stacks \cite{Kresch}
(see specifically section 5.2 there)
has removed the need for this hypothesis.

\section{The boundary of moduli}
\subsection{Boundary of the stack of twisted curves}
We will need to study the geometry of the moduli stack of pre-stable
twisted curves, $\fM^\tw_{g,n}$, as described in
\cite{Olsson-twisted}. In particular we are interested in the
structure of the boundary. We consider the following category
$$\twD(g_1;A\,|\,g_2;B)$$ fibered in groupoids over the category of
schemes. Informally, this category
parametrizes nodal twisted curves, with a distinguished node
separating the curve in two connected components, one of genus $g_1$
containing the markings in a subset $A \subset \{1,\ldots n\}$, the
other, of genus $g_2$, containing the markings in the complementary
set $B$. More formally the objects over a scheme $S$ consist of
commutative diagrams
 $$\xi\ \ \ =\ \ \ \left(
\vcenter{\xymatrix{*+<12pt,12pt>{\cG_1}\ar^{\alpha}[rr] \ar@{^{(}->}[d]& &*+<12pt,12pt>{\cG_2} \ar@{^{(}->}[d] \\ \cC_1 \ar[dr] && \cC_2\ar[dl] \\ &S}}\right)$$
where
\begin{enumerate}
\item $\cC_1 \to S$ is a pre-stable twisted curve of genus $g_1$ with marking in $A\sqcup \checkbullet$,
\item $\cC_2 \to S$ is a pre-stable twisted curve of genus $g_2$ with marking in $B\sqcup \bullet$,
\item $\cG_1$ and $\cG_2$ are the markings on $\cC_1$ and $\cC_2$ corresponding to $\checkbullet$ and $\bullet$, respectively,
and
\item $\alpha$ is an isomorphism inverting the band.
\end{enumerate}
An arrow of $\twD(g_1;A\,|\,g_2;B)$ is a fiber diagram
$$
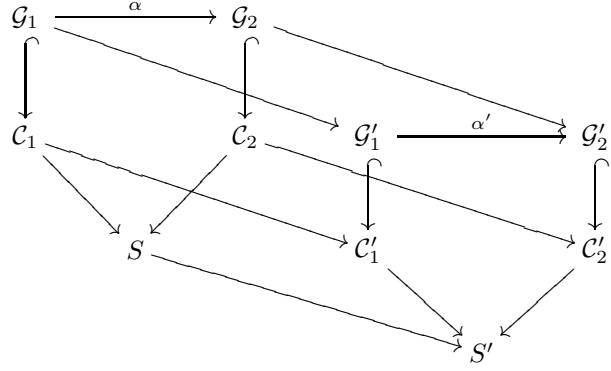
\vcenter{\xymatrix{
*+<12pt,12pt>{\cG_1}\ar^{\alpha}[rr]\ar[rrrd] \ar@{^{(}->}[d]& &*+<12pt,12pt>{\cG_2} \ar@{^{(}->}[d]\ar[rrrd] \\
\cC_1 \ar[dr]\ar[rrrd] && \cC_2\ar[dl]\ar[rrrd] &
*+<12pt,12pt>{\cG'_1}\ar^{\alpha'}[rr] \ar@{^{(}->}[d]& &*+<12pt,12pt>{\cG'_2} \ar@{^{(}->}[d]
\\
&S\ar[rrrd] & &\cC'_1 \ar[dr] && \cC'_2\ar[dl] \\
&&&&S'
}}$$
This, in particular, includes the data of a 2-isomorphism in the following square:
$$\xymatrix{
\cG_1\ar^{\alpha}[r] \ar[d]& \cG_2 \ar[d] \ar@{=>}[dl]\\
 \cG'_1\ar^{\alpha'}[r] & \cG'_2
 }$$
\begin{remark}
In any reasonable framework of 2-stacks, $\twD(g_1;A\,|\,g_2;B)$ should be the fibered product
$$\fM^\tw_{g_1,A\sqcup\checkbullet}\ \ \mathop\times\limits_{\mathop\sqcup\limits_r\fB\cB\bmu_r}\ \ \fM^\tw_{g_2,B\sqcup\bullet}.$$
Here $\fB\cB\bmu_r$ is the classifying 2-stack of $\cB\bmu_r$, which parametrizes gerbes banded by $\bmu_r$,  and the the first morphism underlying the product has the band inverted. Since this 2-stack occurs only as the basis of the fibered product, the result is still a 1-category.

There is a similar construction for non-separating nodes, which we will not describe explicitly here.
\end{remark}

\begin{proposition}
The category $\twD(g_1;A\,|\,g_2;B)$ is a smooth algebraic stack, locally of finite type over the base field $k$.
\end{proposition}

\begin{proof}
Consider the universal gerbes $\cG_1$ and $\cG_2$
corresponding to the markings $\checkbullet$ and $\bullet$ over the
product
$\fM^\tw_{g_1,A\sqcup\checkbullet}\times\fM^\tw_{g_2,B\sqcup\bullet}$.
By Theorem 1.1 of \cite{Olsson-Hom}, there exists an algebraic stack
$$ \Isom_{\fM^\tw_{g_1,A\sqcup\checkbullet}\times\fM^\tw_{g_2,B\sqcup\bullet}}(\cG_1, \cG_2)$$
parametrizing isomorphisms between $\cG_1$ and $\cG_2$. The ``change
of band" isomorphism $\bmu_r\to \bmu_r$  induced by such an
isomorphism $\cG_1 \to \cG_2$ is locally constant, and
$\twD(g_1;A\,|\,g_2;B)$ is the locus where it is the inversion
isomorphism. Since the $\cG_i$ are \'etale gerbes, the smoothness of
the Isom stack follows immediately from that of $\fM^\tw_{g,n}$.
\end{proof}

\begin{proposition}
We have a natural representable morphism
$$gl: \twD(g_1;A\,|\,g_2;B) \to \fM^\tw_{g_1+g_2, A\sqcup B}$$ induced by gluing the
two families of curves over $\twD$ into a family of reducible curves
with a distinguished node.
\end{proposition}
{\bf  Proof.} Fix an object of $\twD(g_1;A\,|\,g_2;B)$ over $S$. By
Proposition \ref{Prop:gluing} applied to the  diagram $$\xymatrix
{*+<12pt,12pt>{\cG_1} \ar@{^{(}->}[r]\ar@{^{(}->}[d] & \cC_2\\\cC_1},$$ we have an
associated family of nodal curves $\cC:= \cC_1 \cup_{\cG_1} \cC_2.$
Representability follows from a straightforward comparison of
isotropy groups. \qed

\begin{definition} \label{Def:gothic-r}
We define the locally constant function
   \[
   \fr : \twD(g_1;A\,|\,g_2;B) \to \ZZ
   \]
(this is a Gothic ``r'') which takes a nodal twisted curve to the index of the node.
\end{definition}

\subsection{Boundary of the stack of twisted stable maps}
There is an analogous gluing map on the spaces of morphisms.

\begin{proposition}\hfill
\begin{enumerate}
\item
Consider the evaluation morphisms
$$ \check e_{ \checkbullet}: \cK_{g_1,A\sqcup \checkbullet}(\cX,\beta_1) \to \riX$$
and
$$            e_{ \bullet}: \cK_{g_2,B\sqcup \bullet}(\cX,\beta_2) \to \riX$$
There exists a natural representable morphism
$$\cK_{g_1,A\sqcup \checkbullet}(\cX,\beta_1) \times_{\riX}\cK_{g_2,B\sqcup \bullet}(\cX,\beta_2)  \ \ \lrar \ \  \cK_{g_1+g_2,A\sqcup B}(\cX,\beta_1+\beta_2).$$
\item
Consider the evaluation morphisms
$$ \check e_{\checkbullet}\times  e_{\bullet}: \cK_{g-1,A\sqcup\{\checkbullet,\bullet\}}(\cX,\beta) \to \riXs$$
There exists a natural representable morphism
$$\cK_{g-1,A\sqcup\{\checkbullet,\bullet\}}(\cX,\beta) \times_{\riXs}\riX  \ \ \lrar \ \  \cK_{g,A}(\cX,\beta).$$
\end{enumerate}
\end{proposition}

\begin{proof}
We prove the first statement, the second being similar, replacing Proposition \ref{Prop:gluing} with Corollary \ref{Cor:limits}.

We give the morphism on the level of $S$-valued points. We have an identification of objects
\begin{equation}\label{Diag:esuoh}
\cK_{g_1,A\sqcup \checkbullet}(\cX,\beta_1) \mathop\times\limits_{\ \riX}\cK_{g_2,B\sqcup \bullet}(\cX,\beta_2)\ \big(\  S \ \big) \ \ = \ \ \left\{
\vcenter{\xymatrix{*+<12pt,12pt>{\cG_1}\ar^{\alpha}[rr] \ar@{^{(}->}[d]& &*+<12pt,12pt>{\cG_2} \ar@{^{(}->}[d] \\ \cC_1 \ar[dr]\ar[r] &\cX& \cC_2\ar[l]\ar[dl] \\ &S}}\right\}
\end{equation}
where the diagram is a 1-commutative diagram of stacks, $\cG_i \subset \cC_i$ are the markings corresponding to $\checkbullet$ and $\bullet$, respectively, and $\alpha: \cG_1 \to \cG_2$ is the isomorphism \emph{inverting the band} induced by $ \check e_{\checkbullet}$ and $ e_{\bullet}$.

By Proposition \ref{Prop:gluing} applied to the  diagram $$\xymatrix {*+<12pt,12pt>{\cG_1} \ar@{^{(}->}[r]\ar@{^{(}->}[d] & \cC_2\\\cC_1},$$ we have an associated
family of nodal curves $\cC:= \cC_1 \cup_{\cG_1} \cC_2.$ Since Diagram
(\ref{Diag:esuoh}) is commutative, we are given a 2-isomorphism between the two resulting maps $\cG_1 \to \cX$.  Therefore, by the universal property of $\cC$, we have a morphism $\cC \to \cX$, which is clearly a twisted stable map over $S$.
\end{proof}

\begin{proposition}\label{cartesian}
We have a cartesian diagram
$$\xymatrix{{\mathop\coprod\limits_{\beta_1+\beta_2 = \beta} \cK_{g_1,A\sqcup \checkbullet}(\cX,\beta_1) \times_{\riX}\cK_{g_2,B\sqcup \bullet}(\cX,\beta_2) \ar[r]\ar[d]} &
\cK_{g_1+g_2,A\sqcup B}(\cX,\beta)\ar[d] \\
\twD(g_1;A\,|\,g_2;B)\ar^{gl}[r] &\fM^\tw_{g_1+g_2, A\sqcup B}
}$$
\end{proposition}

\begin{proof}
For convenience of notation we will use the shorthand
$$\xymatrix{{ \mathop\coprod\limits_{\beta_1+\beta_2 = \beta}\cK_1(\beta_1)\times_{\riX}\cK_2(\beta_2) \ar[r]\ar[d]} &
\cK(\beta)\ar[d] \\
\twD\ar[r] &\fM^\tw
}$$

The diagram gives a morphism $$\mathop\coprod\limits_{\beta_1+\beta_2 = \beta} \cK_1(\beta_1) \times_{\riX}\cK_2(\beta_2)  \to \cK(\beta) \times_{\fM^\tw} \twD. $$ We construct a morphism in the reverse direction as follows. It suffices to restrict attention to points over a connected base scheme $S$. In this case an object on the right hand side is a triple
$$ \left( f:\cC \to \cX,\ \xi, \ \phi\right) ,$$
where $\xi$ is an object of $\twD(S)$ and $\phi$ is an isomorphism between the resulting objects in $\twM$, namely between $\cC$ and $\cC_1\cup_{\cG_1}\cC_2$. Since $f:\cC \to \cX$ is stable, so are the resulting morphisms $\cC_i\to \cX$, with the additional markings taken into account. By the connectedness of the base, these maps have  constant image classes $\beta_1,\beta_2$.   With this, the diagram describing $\xi$ is completed to a diagram as in the description \ref{Diag:esuoh} of an object on the left hand side, which is what we needed.
\end{proof}

\subsection{Gluing and virtual fundamental classes}
For Gromov-Witten theory, one of the key points is that a similar
statement holds for the virtual fundamental class. First, a crucial
fact is that fundamental classes exist for the stacks of twisted
curves we need. This is because, as mentioned in Section
\ref{Sec:stack-twisted-curves}, away from a closed substack of
arbitrarily high codimension, the stack of twisted curves of bounded
indices is of finite type.  We will show below (Lemma
\ref{Lem:gluing-finite}) that  $gl:\twD \to \twM$ is a finite
unramified  morphism. By \cite{Kresch}, Section 4.1 it induces a
pull-back homomorphism on Chow groups
$$gl^! : A_*(\cK(\beta)) \lrar \mathop\oplus\limits_{\beta_1+\beta_2=\beta}
A_*(\cK_1(\beta_1)\times_{\riX} \cK_2(\beta_2)).$$  Pulling back the
virtual class on $\cK(\beta)$ gives us a candidate for the virtual
fundamental class of the boundary. There is another natural Chow
class living on the fibered product spaces coming from the pull back
by the diagonal morphism $\Delta: \riX \to \riXs$. The splitting
axiom in Gromov-Witten theory identifies these two classes.

\begin{proposition}\label{Prop:virtual}
$$gl^![\cK_{g,A\cup B}(\cX,\beta)]^\vir = \sum_{\beta_1+\beta_2=\beta}
\Delta^!([\cK_{g_1,A\sqcup \checkbullet}(\cX,\beta_1)]^\vir \times
[\cK_{g_2,B\sqcup \bullet}(\cX,\beta_2)]^\vir).$$
\end{proposition}

\begin{proof}
This proof is very similar to the proof of the splitting axiom
for Gromov-Witten theory of schemes.  First, we observe that the
left hand side of our equation is the
pullback of a relative virtual fundamental class under a change of base.
It follows by Proposition 7.2 of \cite{BF}
that this left hand side is the relative virtual fundamental class
of $\cK_1 \times_{\riX} \cK_2$ over $\twD $ with respect to the
relative perfect obstruction theory $R\pi_*(f^*T\cX)$.

We need to compare this to the right hand side. Here we use the
basic compatibility result for virtual fundamental classes. The
class $[\cK_1]^\vir \times [\cK_2]^\vir$ is the virtual fundamental
class associated with the relative perfect obstruction theory on
$\cK_1 \times \cK_2$ given by $R\pi_{1*}(f_1^*T\cX)\oplus
R\pi_{2*}(f_2^*T\cX)$.  This is an immediate consequence of
Proposition 5.7 of \cite{BF}. By considering the normalization
sequence for the family of nodal curves with a distinguished node
over $\cK_1 \times_{\riX} \cK_2$, we get the following distinguished
triangle:
$$ R\pi_*(f^*T\cX) \to R\pi_{1*}(f_1^*T\cX)\oplus R\pi_{2*}(f_2^*T\cX)
\to R\pi_{\Sigma *} (f_\Sigma^*T\cX) $$ where $f_\Sigma$ denotes the
restriction of $f$ to the gerbe which is the intersection of $\cC_1$
with $\cC_2$. By Proposition 5.10 of \cite{BF}, in order to prove
the equality we want, we just need to identify $R\pi_{\Sigma
*}(f_\Sigma^* T\cX)$ with the normal bundle of the map $\Delta$. The
normal bundle of $\Delta$ is obviously $T\riX$. Applying Lemma
\ref{tangent} to $S=\cK_1\times_{\riX}\cK_2$ with the morphism
$F=f_\Sigma$ gives our result.
\end{proof}

An identical argument yields the analogous splitting axiom for a
non-separating node.
\begin{proposition}
$$gl^![\cK_{g,A}(\cX,\beta)]^\vir =
\Delta^![\cK_{g-1,A\sqcup \{ \bullet, \checkbullet\}
}(\cX,\beta_1)]^\vir.$$
\end{proposition}

\section{Gromov--Witten classes}

\subsection{Algebraic Gromov--Witten classes}

\begin{definition}
We define a locally constant function
$r: \riX \to \ZZ$ by evaluating on geometric points:
$(x,\cG) \mapsto r$,  where $\cG$ is a gerbe banded by $\bmu_r$.  We view $r$ as an
element in $A^0(\riX)$.
\end{definition}

We now define \emph{Gromov--Witten Chow classes:}

\begin{definition}\label{Def:gromov-witten-classes} Fix integers $g,n$, Chow classes
$\gamma_i\in A^*(\riX)_\QQ$, and a curve class $\beta$. We define a
class in $A_*(\riX)_\QQ$ by the formula
$$\langle \gamma_1,\ldots,\gamma_n, * \rangle_{g,\beta}^\cX \ \ = \ \
r\ \cdot\ \check e_{n+1\,*}\left( \left(\prod_{i=1}^n e_i^*\gamma_i\right)\cap\left[\cK_{g,n+1}(\cX,\beta)\right]^\vir\right).$$
\end{definition}
We will suppress the superscript $\cX$ when the target is clear, and the genus $g$ when $g=0$.
We will write $\langle \gamma_1, \gamma_2,\delta_A, * \rangle_{g,\beta}^\cX$
for an expression of the type $\langle \gamma_1,
\gamma_2,\delta_{i_1},\ldots,\delta_{i_m}, * \rangle_{g,\beta}^\cX$
with $A = \{i_1,\ldots i_m\}$ subject to the convention $i_1<\cdots<
i_m$. The factor $r$ comes very
naturally in the proof, though see Proposition
\ref{Prop:lifted-evaluation} below for a way to avoid this
factor.

\subsubsection{Alternative formalism} \label{Sec:formalism}
There is at least one other way to define Gromov\ddash Witten classes
 introduced in \cite{AGV}, and it is necessary to compare
 them. In
 that paper we defined $$\ocM_{g,n}(\cX,\beta) = \Sigma_1^\cC
 \mathop\times\limits_{\cK_{g,n}(\cX,\beta)} \cdots
 \mathop\times\limits_{\cK_{g,n}(\cX,\beta)} \Sigma_n^\cC,$$ the
 moduli stack of twisted stable maps with sections at the markings. It
 has degree $(\prod r_i)^{-1}$ over  $\cK_{g,n}(\cX,\beta)$. There are
 clearly natural evaluation maps 
$$
{e}^\cM_i:\ocM_{g,n}(\cX,\beta) \to \iX
$$
to the \emph{non-rigidified}
 cyclotomic inertia stack.
Given $n$ Chow
 classes  $\widetilde\gamma_i\in A^*(\iX)_\QQ$ on the non-rigidified
 cyclotomic inertia stack, we defined classes
$$\langle \widetilde\gamma_1,\ldots,\widetilde\gamma_n, *
 \rangle_{g,\beta}^\cX \ \ = \ \
 {\check{{e}}}^\cM_{n+1\,*}\left( (\prod r_i) \left(\prod_{i=1}^n
\tilde{e}^{\cM\,*}_i\widetilde\gamma_i\right)\cap
 \left[\ocM_{g,n+1}(\cX,\beta)\right]^\vir\right).$$  This formalism is
 used in the work of Tseng \cite{Tseng}.

There is another way to write down the same classes without need to
introduce the stack
$\ocM_{g,n}(\cX,\beta)$. Even though liftings $\tilde{e}_i:
\cK_{g,n}(\cX,\beta) \to \iX$ do not necessarily exist, $$\xymatrix{ &
  \iX\ar[d]^{\varpi}\\
\cK_{g,n}(\cX,\beta)\ar[r]_{e_i}\ar@{.>}[ru]|{\not\exists\, \tilde e_i} &  \riX
}$$
the
isomorphism between the rational Chow groups (or cohomology groups) of
$\iX$ and $\riX$ enables us to move from one to the other on the
intersection theory level.
A  lifting $\tilde e_{i\,*}$ of $e_{i\,*}$ is obtained by composing
with the non-multiplicative  isomorphism
$$(\varpi_*)^{-1}: A^*(\riX)_\QQ \to A^*(\iX)_\QQ.$$ Note
$(\varpi_*)^{-1} = r \cdot \varpi^*$, so $$\tilde e_{i\,*} =
(\varpi_*)^{-1} \circ e_{i\,*} = r \varpi^*\circ e_{i\,*}.$$
Similarly we define  $$\tilde e_i^* = e_i^*\circ (\varpi^*)^{-1}.$$
Since $(\varpi^*)^{-1}= r\cdot \varpi_*$ we can also write $\tilde
e_i^* = r\cdot e_i^*\circ \varpi_*.$  (We remark that the
corresponding formula in \cite{AGV} has $r$ mistakenly replaced by
$r^{-1}$.)

The basic comparison result is

\begin{proposition}\label {Prop:lifted-evaluation}
\begin{enumerate}
\item  For Chow classes $\gamma_i\in A^*(\riX)_\QQ$, we have
$$\varpi^*\langle \gamma_1,\ldots,\gamma_n, *
  \rangle_{g,\beta}^\cX \ \ = \ \ \langle
  \varpi^*\gamma_1,\ldots,\varpi^*\gamma_n, *
  \rangle_{g,\beta}^\cX $$
\item For Chow classes $\widetilde\gamma_i\in A^*(\iX)_\QQ$, we have
$$\langle \widetilde\gamma_1,\ldots,\widetilde\gamma_n, *
  \rangle_{g,\beta}^\cX \ \ = \ \
 {{(\iota\circ \tilde{e}_{n+1})}}_{*}\left( \left(\prod_{i=1}^n
\tilde{e}_i^*\widetilde\gamma_i\right)\cap
 \left[\cK_{g,n+1}(\cX,\beta)\right]^\vir\right).$$
\end{enumerate}
\end{proposition}

The first part shows that, if one  identifies $A^*(\riX)_\QQ$ and
$A^*(\iX)_\QQ$ \emph{using the multiplicative homomorphism $\varpi^*$},
the Gromov--Witten classes are unchanged. The second part says that,
if one is willing to pretend a lifting $\tilde
e_i:\cK_{g,n}(\cX,\beta) \to \iX$ exists, all the factors of $r$ are
removed from the formalism. While the rigidified inertia stack and  evaluation map  $
e_i:\cK_{g,n}(\cX,\beta) \to \riX$ is what arises naturally,  the
formalism using $\tilde
e_i$ and $\iX$  is probably the most convenient
one to work with, and has been used  in the work of Cadman
\cite{Cadman1},\cite{Cadman2}. We will use this formalism in our
example in section \ref{Sec:example}. A direct comparison was carried
out in  the example in \cite{cimenotes}.

\begin{proof}
This is immediate using the non-cartesian commutative diagram
   \[
   \xymatrix{
   \save+<0pt,+10pt>*+{\cM_{g,n+1}(\cX,\beta)}\ar_\rho[d]\ar[r]^(0.6){\tilde e^\cM_i}\restore
    &\iX\ar^\varpi[d]\\ \cK_{g,n+1}(\cX,\beta)\ar[r]^(0.6){e_i} \ar@{.>}[ru]^{\tilde e_i}
    & \riX,
   }
   \]
where $\deg \rho = (\prod r_i)^{-1}$ and   $\deg \varpi = 1/r_i$.
\end{proof}

\subsection{The WDVV equation}

In genus 0 we have the \emph{Witten--Dijkgraaf--Verlinde--Verlinde} (WDVV) equation:
\begin{theorem}\label{Thm:WDVV}
\begin{eqnarray*}
\lefteqn{\sum_{\beta_1+\beta_2 = \beta}\ \  \sum_{A \sqcup B = \{1,\ldots ,n\}}
\left\langle \left\langle \gamma_1, \gamma_2,\delta_A,*\right\rangle_{\beta_1},\gamma_3, \delta_B,*\right\rangle_{\beta_2}  = } \\
& &\sum_{\beta_1+\beta_2 = \beta}\ \  \sum_{ A \sqcup B = \{1,\ldots ,n\}}
\left\langle \left\langle \gamma_1, \gamma_3,\delta_A,*\right\rangle_{\beta_1},\gamma_2, \delta_B,*\right\rangle_{\beta_2}
\end{eqnarray*}
\end{theorem}

Note that the two lines differ precisely in the positions of $\gamma_2$ and $\gamma_3$.

\begin{proof}
Consider the stabilization morphism
$$st: \cK_{0,n+4}(\cX,\beta) \to \ocM_{0,4}$$ corresponding to forgetting the map to $\cX$, passing to coarse curves, and forgetting  all the marking except the last four.

We show the equality by showing that both sides equal
$$\Psi \ = \ r \ \cdot \ \check e_{n+4\ *}\left( \left (st^*[pt] \cup \prod_{i=1}^n e_i^*\delta_i \cup \prod_{j=1}^3 e_{n+j}^*\gamma_j\right)\cap\left[\cK_{0,n+4}(\cX,\beta)\right]^\vir\right).$$

Write $\hat A = A \sqcup\{n+1,n+2\}$ and $\hat B = B \sqcup\{n+3,n+4\}$
Consider the following  diagram:
$$\xymatrix{\coprod \cK_{0,\hat A\sqcup \checkbullet}(\cX, \beta_1) \times_{\riX}\cK_{0,\hat B \sqcup \bullet}(\cX, \beta_2) \ar^(.8)l[r] \ar[d]
    & \cK\ar[d]\ar@/^2pc/[dd]^{st}\\
\coprod \twD(\hat A|\hat B)\ar^{gl}[r]& \fM^\tw_{0,n+4}\ar[d]
\\
&  \ocM_{0,4}
}$$

\begin{proposition}\label{Prop:divisor}
$$st^*[pt] \cap [\cK]^\vir = l_*\left( e^*_\bullet r \ \cdot \ gl^![\cK]^\vir\right).$$
\end{proposition}
\begin{proof}
We expand the diagram into the following cartesian diagram:
$$\xymatrix{
{}\save+<0pt,-0pt>
*+{\coprod^{\vphantom{\riX}} \cK_{0,\hat A \sqcup \checkbullet}(\cX, \beta_1) \mathop\times\limits_{\riX}\cK_{0,\hat B \sqcup \bullet}(\cX, \beta_2)}
\ar^l[rrr] \ar_\phi[d]
\restore
   & &  & \cK\ar[d]\ar@/^2pc/[dddd]^{st}\\
\coprod \twD(\hat A|\hat B)\ar^{gl}[rrr]\ar[d] &&& \fM^\tw_{0,n+4}\ar[d]
\\
\twD(12|34)\ar^{i_\fP}[r] &\fP\ar[d]\ar[rr] && \fM^\tw_{0,4}\ar[d]
 \\
& \fM_{0,3}\times \fM_{0,3} \ar^(.7){i_\fQ}[r] & \fQ \ar[r]\ar[d] & \fM_{0,4} \ar[d]
  \\
  && \{pt\}\ar[r] &     \ocM_{0,4}
}$$
\begin{lemma} All the morphisms in
$$\fM^\tw_{0,n+4} \lrar \fM^\tw_{0,4}\lrar \fM_{0,4} \lrar \ocM_{0,4}
$$
are flat.
\end{lemma}

\begin{proof}
The first arrow just forgets the first $n$ markings which is smooth.  The second is locally a finite morphism of smooth stacks. The third is a dominant morphism from a smooth stack to a smooth curve.
\end{proof}

\begin{lemma}\label{Lem:gluing-finite}
All horizontal arrows in this diagram are finite and unramified.
\end{lemma}
\begin{proof}
We first show the result for the arrow
$\fM_{0,3}\times \fM_{0,3} \to\fM_{0,4}$. Consider the universal
curve $f_{0,4}:\fC_{0,4} \to \fM_{0,4}$, which is clearly a proper
and representable morphism. There is a closed substack
$\Sing(f_{0,4}) \subset \fC_{0,4}$ consisting of the nodes of the
universal curve, schematically defined by the first Fitting ideal of
$\Omega^{1}_{\fC_{0,4}/\fM_{0,4}}$. Since $\Sing(f_{0,4}) \to
\fM_{0,4}$ is representable, quasi-finite and proper, it is a finite
morphism.

Inside $\Sing(f_{0,4})$ we have an open and closed substack
$\Sigma\subset \Sing(f_{0,4})$ consisting of nodes separating the
marking numbered $1,2$ from the markings numbered $3,4$. We claim
that there is an isomorphism $\Sigma\simeq \fM_{0,3}\times
\fM_{0,3}$ over $\fM_{0,4}$, which proves the required property.

We construct a morphism $\Sigma \to \fM_{0,3}\times \fM_{0,3}$ as
follows. The pullback $\fC_{\Sigma} \to \Sigma$ of the universal
pre-stable curve has a canonical section landing at the appropriate
node. The normalization of $\fC_{\Sigma}$ is the disjoint union of
two families of $3$-pointed curves, obtained by separating the node
marked by the section above. This gives the required morphism.

A morphism in the other direction can be constructed as follows. We
have two universal families $\fC_{1} = \fC_{0,3} \times \fM_{0,3}$
and $\fC_{2} = \fM_{0,3} \times
\fC_{0,3}$. We have a non-cartesian diagram
   \[
   \xymatrix{
   \save+<0pt,+10pt>*+{\fC_{1}\sqcup \fC_{2}}\ar[d]\ar[r]\restore
      &\fC_{0,4}\ar[d]\\
   \fM_{0,3}\times \fM_{0,3}\ar[r] & \fM_{0,4}
   }
   \]
where the diagonal arrow is the gluing map of the third marking of $\fC_{1}$ with the first marking of $\fC_{2}$. Composing the third section $ \fM_{0,3}\times \fM_{0,3} \to \fC_{1}$ with this diagonal arrow (or, for that matter, using the first section of  $\fM_{0,3}\times \fM_{0,3} \to \fC_{2}$), we get a morphism $ \fM_{0,3}\times \fM_{0,3} \to \fC_{0,4}$, which obviously lands in $\Sigma$. It is a simple local computation to check that the two arrows are inverses of each other.

 The arrow $i_\fP$ is the embedding of the reduced substack in $\fP$.  All other arrows
 arise by base change.
\end{proof}

 As a consequence, we can make sense of the following lemma.
 \begin{lemma}
$$i_{\fQ\, *} \left[ \fM_{0,3}\times \fM_{0,3}\right] \ \ = \ \ [\fQ]$$
and, using Definition \ref{Def:gothic-r},
$$i_{\fP\, *}\left( \fr \ \cdot    \left[ \twD(12|34)\right]\right) \ \ = \ \ [\fP]$$
\end{lemma}

 The first statement is well-known---see \cite{Behrend}, Proposition 8.  The second follows from the deformation theory of twisted curves, as mentioned in section \ref{Sec:stack-twisted-curves}. Proposition \ref{Prop:divisor} now follows from the lemmas using the projection formula and the observation that $e^*_\bullet r = e^*_{\checkbullet} r = \phi^*\fr$.
\end{proof}

\begin{corollary}
\begin{eqnarray*}
\Psi &=& \ r \ \cdot \ \check e_{n+4\ *}\left( \prod_{i=1}^n e_i^*\delta_i \cup \prod_{j=1}^3 e_{n+j}^*\gamma_j\right)\cap l_*\left( e^*_\bullet r \ \cdot \ gl^![\cK]^\vir\right)\\
&=&\ r \ \cdot \ (\check e_{n+4}\circ l)_*\ \left( e^*_\bullet r \ \cdot \prod_{i=1}^n l^*e_i^*\delta_i \cup \prod_{j=1}^3 l^*e_{n+j}^*\gamma_j\ \cap \ gl^![\cK]^\vir\right)
\end{eqnarray*}
\end{corollary}

To simplify notation, we fix $$\cK_1= \cK_{0,\hat A \sqcup \checkbullet}(\cX, \beta_1) , \ \cK_2= \cK_{0,\hat B \sqcup \bullet}(\cX, \beta_2) ,
\text{   and   } \  \cK= \cK_{0,n+4}(\cX, \beta) $$
and $$\eta_1= e_{n+1}^*\gamma_1\cup e_{n+2}^*\gamma_2\cup \prod_{i\in A}e_i^*\delta_i \ \ \text{   and   } \ \
\eta_2 = e_{n+3}^*\gamma_3\cup \prod_{i\in B} e_i^*\delta_i .$$
Consider the following diagram.
$$\xymatrix{\cK_1 \times_{\riX} \cK_2 \ar[r]^(.6){p_2}\ar[d]_{p_1}& \cK_2 \ar[r]^{\check e_{n+4}}\ar[d]^{e_\bullet} & \riX\\
\cK_1\ar[r]^{\check e_{\checkbullet}} & \riX}$$

The expression $\left\langle \left\langle \gamma_1, \gamma_2,\delta_A,*\right\rangle_{\beta_1},\gamma_3, \delta_B,*\right\rangle_{\beta_2}$ is by definition
$$r\ \cdot \ \check e_{n+4\, *}\left(\eta_2\cup e_\bullet ^*\left(r\cdot \check e_{\checkbullet\, *}\left(\eta_1\cap[\cK_1]^\vir\right)\right)\cap [\cK_2]^\vir\right)
.$$
By the following lemma, this expression equals
$$r \cdot \check e_{n+4 \,*}(e_\bullet^*r \ \cdot\  p_{2\,*}(p_2^*\eta_2 \cup p_1^*\eta_1 \cap \Delta^! ([\cK_1]^\vir\times[\cK_2]^\vir)).$$
Applying Proposition \ref{Prop:virtual} and summing over $A,B$ and $\beta_1,\beta_2$ gives the Theorem.
\end{proof}

\begin{lemma}\label{Lem:virtual-gysin}
   \begin{align*}
      \eta_2\cup & e_\bullet ^*\left(\check e_{\checkbullet\, *}\left(\eta_1\cap[\cK_1]^\vir\right)\right)\cap [\cK_2]^\vir\\
   &= p_{2\,*}(p_2^*\eta_2 \cup p_1^*\eta_1 \cap \Delta^! ([\cK_1]^\vir\times[\cK_2]^\vir).
   \end{align*}
\end{lemma}

\begin{proof}
Set $\xi_{i} = \eta_{i}\cap [\cK_{i}]^{\vir}$; then the left hand side of the equality is
   \[
   e_\bullet ^*\left(\check e_{\checkbullet\, *}\left(\xi_{1}\right)\right)\cap \xi_{2},
   \]
while the right hand side is
   \[
   p_{2\,*}\Delta^! (\xi_{1}\times\xi_{2}).
   \]
Consider the following cartesian diagram
   \[
   \xymatrix{
   \cK_{1}\times_{\riX}\cK_{2} \ar[r]^{\qquad p_{2}}\ar[d] & \cK_{2} \ar[r]^{e_{\bullet}} \ar[d]^{\Gamma_{e_{\bullet}}} & \riX \ar[d]^{\Delta}\\
   \cK_{1}\times\cK_{2} \ar[r]^{\check e_{\checkbullet}\times \mathrm{id}\quad} \ar[d]_{\pi_{\cK_{1}}} & \riX\times \cK_{2} \ar[r]^{\mathrm{id} \times e_{\bullet}\quad}\ar[d]^{\pi_{\riX}} & \riX\times \riX \\
   \cK_{1}\ar[r]^{\check e_{\checkbullet}} &\riX
   }
   \]
We have
  \begin{align*}
  e_\bullet ^*\left(\check e_{\checkbullet\, *}\left(\xi_{1}\right)\right)\cap \xi_{2}
&= \Gamma_{e_\bullet}^* \left( \pi_{\riX}^{*}\check e_{\checkbullet\, *}\xi_{1}\right)\ \cap \ \xi_{2} \\
&= \Gamma_{e_\bullet}^* \left(\check e_{\checkbullet\, *}\xi_{1}\ \times \ \xi_{2}\right) \\
&=\Gamma_{e_\bullet}^* \left( (\check e_{\checkbullet}\times \mathrm{id})_{*}(\xi_{1}\ \times \ \xi_{2})\right) \\
&=p_{2\,*}\Gamma_{e_\bullet}^! \left(\xi_{1}\ \times \ \xi_{2}\right) \\
&=p_{2\,*}\Delta^{!} \left(\xi_{1}\ \times \ \xi_{2}\right).\qedhere
  \end{align*}
\end{proof}

\subsection{Topological Gromov--Witten classes}

In this section the base field will be $k = \CC$.
The definition of a Gromov--Witten cohomology class is the same as in Definition~\ref{Def:gromov-witten-classes}:

\begin{definition} Fix integers $g,n$, classes $\gamma_i\in \H^*(\riX)_\QQ$, and a curve class $\beta\in N^{+}(X)$. We define
$$\langle \gamma_1,\ldots,\gamma_n, * \rangle_{g,\beta}^\cX \ \ = \ \
r\ \cdot\ \check e_{n+1\,*}\left( \left(\prod_{i=1}^n e_i^*\gamma_i\right)\cap\left[\cK_{g,n+1}(\cX,\beta)\right]^\vir\right)$$
where by $\left[\cK_{g,n+1}(\cX,\beta)\right]^\vir$ we mean the homology class in $\H_{*}\bigl(K_{g,n+1}(\cX,\beta),\QQ\bigr)$ corresponding to the virtual fundamental class
   \[
   \left[\cK_{g,n+1}(\cX,\beta)\right]^\vir \in
      A_{*}\bigl(K_{g,n+1}(\cX,\beta)\bigr)_\QQ.
   \]
\end{definition}

In genus 0 we have again the WDVV equation. To state it precisely we
need a sign convention.

We restrict the discussion to cohomology classes
$\gamma_i$ and $\delta_i$ which are homogeneous with respect to the
usual grading in $\H^*(\riX)$ - in fact homogeneous parity
suffices; the formulas extend to the inhomogeneous case, but are less clean.
For $A\sqcup B = \{1,\ldots ,n\}$  we can write $\delta_A =
\delta_{i_1}\wedge \cdots \wedge \delta_{i_m}$, where
$A = \{i_1,\ldots i_m\}$ subject to the ordering convention $i_1<\cdots<
i_m$, and similarly for $\delta_B$.   We define signs
$(-1)^{\epsilon_1(A)}$ and  $(-1)^{\epsilon_2(A)}$  so that
$$
\left( \gamma_1 \wedge  \gamma_2\wedge
 \gamma_3\right)\ \wedge\ \left(\delta_1 \wedge \cdots\wedge
 \delta_n\right)\ \ =\ \ (-1)^{\epsilon_1(A)} \left( \gamma_1 \wedge
 \gamma_2 \wedge \delta_A\right)\ \wedge\
\left( \gamma_3 \wedge \delta_B\right)
$$
and
$$
\left( \gamma_1 \wedge  \gamma_2\wedge
 \gamma_3\right)\ \wedge\ \left( \delta_1 \wedge \cdots\wedge
 \delta_n\right)\ \ =\ \ (-1)^{\epsilon_2(A)} \left( \gamma_1 \wedge
 \gamma_3 \wedge \delta_A\right)\ \wedge\
\left( \gamma_2 \wedge \delta_B\right),
$$
Of course the products could vanish, but the signs are formally well defined depending only on the parity of the classes $\gamma_i$ and $\delta_i$.
With these conventions we have:

\begin{theorem}\label{Thm:WDVV-homology}
\begin{eqnarray*}
\lefteqn{\sum_{\beta_1+\beta_2 = \beta}\ \  \sum_{A \sqcup B = \{1,\ldots ,n\}}
(-1)^{\epsilon_1(A)} \left\langle \left\langle \gamma_1, \gamma_2,\delta_A,*\right\rangle_{\beta_1},\gamma_3, \delta_B,*\right\rangle_{\beta_2}  = } \\
& &\sum_{\beta_1+\beta_2 = \beta}\ \  \sum_{A \sqcup B = \{1,\ldots
,n\}} (-1)^{\epsilon_2(A)} \left\langle \left\langle \gamma_1,
\gamma_3,\delta_A,*\right\rangle_{\beta_1},\gamma_2,
\delta_B,*\right\rangle_{\beta_2}.
\end{eqnarray*}
\end{theorem}

The proof is identical to that of WDVV in the algebraic case.

We stress that the degrees used in determining the signs above are the
standard degrees in cohomology, which may be different from those in
the age grading defined in the next section.

\subsection{Gromov--Witten numbers}

The usual definition of Gromov--Witten numbers works without changes in our context.

\begin{definition}
 Fix integers $g,n$, classes $\gamma_i\in \H^*(\riX)_\QQ$, and a curve class $\beta\in N^{+}(X)$. We define
$$\langle \gamma_1,\ldots,\gamma_n \rangle_{g,\beta}^\cX \ \ = \ \
\int_{K_{g,n}(\cX,\beta)}\left( \left(\prod_{i=1}^n e_i^*\gamma_i\right)\cap\left[\cK_{g,n}(\cX,\beta)\right]^\vir\right)$$
\end{definition}

The two definitions are connected by the following Proposition.
First, some notation:  let $\alpha_{1},\ldots,\alpha_{M}$ be a basis
for the cohomology of $\riX$. We write $g_{ij} = \int_{\riX}
\alpha_{i}\cap \iota^{*}\alpha_{j}$ and denote by $g^{ij}$ the inverse
matrix.

\begin{proposition}\hfil
\begin{enumerate}

\item
   $\displaystyle
   \langle \gamma_1,\ldots,\gamma_n \rangle_{g,\beta}^\cX =
      \int_{\riX}\frac{1}{r}
         \langle\gamma_1,\ldots,\gamma_{n-1}, * \rangle_{g,\beta}^\cX
         \cap \iota^{*}(\gamma_{n}).
   $
\item
   $\displaystyle
        \langle\gamma_1,\ldots,\gamma_{n-1}, * \rangle_{g,\beta}^\cX\ =\ r\cdot \sum_{i,j=1}^{M} \langle \gamma_1,\ldots,\gamma_{n-1},\alpha_{i} \rangle_{g,\beta}^\cX\  g^{{ij}}\  \alpha_{j}
   $

\end{enumerate}
\end{proposition}

\begin{proof}
This is immediate from the projection formula, noting
that $\check e_{n} = \iota \circ e_{n}$.
\end{proof}

As in Section \ref{Sec:formalism} and Proposition
\ref{Prop:lifted-evaluation} we can again make these formulas
appear even more analogous
to the usual manifold case (i.e. we can remove the factors of $r$)
by multiplying the intersection form on $\H^*(\riX[r])$ by $1/r$. We
will denote this modified intersection form by $\tilde g_{ij}$, and correspondingly define $\tilde g^{ij}=r\cdot g^{ij}$.  This
is equivalent to pulling back the classes to $\iX$ (or considering
directly classes on $\iX$) and doing the
intersection there.

 We can now state the WDVV equation in its
classical form:
\begin{theorem}\label{bigwdvv}
\begin{eqnarray*}
\lefteqn{\sum_{\substack{\beta_1+\beta_2 = \beta \\ A \sqcup B = \{1,\ldots ,n\}}}\sum_{i,j=1}^{M}
 (-1)^{\epsilon_1(A)} \left\langle \delta_A,\gamma_1, \gamma_2,\alpha_{i}\right\rangle_{\beta_1} \ \tilde g^{ij}\ \left\langle\alpha_j,\delta_B,\gamma_3, \gamma_{4}\right\rangle_{\beta_2}  = } \\
& &\sum_{\substack{\beta_1+\beta_2 = \beta\\ A \sqcup B = \{1,\ldots ,n\}}}\sum_{i,j=1}^{M}
 (-1)^{\epsilon_2(A)} \left\langle \delta_A,\gamma_1, \gamma_3,\alpha_{i}\right\rangle_{\beta_1}\ \tilde g^{ij}\ \left\langle\alpha_j, \delta_B,\gamma_2, \gamma_{4}\right\rangle_{\beta_2}.
 \end{eqnarray*}
\end{theorem}
\section{The age grading}

\subsection{The age of a sheaf}

Consider the group-scheme $\bmu_{r}$ over a field, with its
representation ring $\mathrm{R}\bmu_{r}$. Each character
$\lambda\colon \bmu_{r} \to \bG_{\mathrm{m}}$ is of the form $t
\mapsto t^{k}$ for a unique integer $k$ with $0 \leq k \leq r-1$;
following M. Reid (see e.g. \cite{Ito-Reid}), we define the \emph{age} of
$\lambda$ as $k/r$. Since these characters form a basis for the
representation ring of $\bmu_{r}$, this extends to a unique additive
homomorphism $\age \colon \mathrm{R} \bmu_{r} \to \QQ$.

Now let $\cG \to T$ be a gerbe banded by $\bmu_{r}$, and let $\cE$ be a locally free sheaf on $\cG$. There is an \'etale covering $\{T_{i}\to T\}$ with sections $T_{i} \to \cG$, inducing an isomorphism $\cG_{T_{i}} \simeq \cB(\bmu_{r})_{T_i}$. Then the pullback of $\cE$ to $\cG_{T_{i}}$ becomes a locally free  sheaf $\cE_{T_{i}}$ on $T_{i}$ with an action of $\bmu_{r}$; and the age of each fiber is a locally constant invariant. Furthermore, this invariant is independent of the section, and the age of $\cE$ is a locally constant function on $T$.

Consider a connected scheme $T$, and an object of $\riX(T)$, consisting of a gerbe $\cG \to T$ and a representable morphism $f:\cG \to \cX$. Then the \emph{age} of the object is a rational number defined to be the age of the locally free sheaf $f^{*}\mathrm{T}_{\cX}$. This number only depends on the connected component of $\riX$ containing the image of $T$.

\begin{definition}
The age of a connected component $\Omega$ of $\riX$ is the age of any object of $\Omega(T)$, where $T$ is a connected scheme.
\end{definition}

The age is called the \emph{degree-shifting number} in \cite{Chen-Ruan}.

\subsection{Riemann-Roch for twisted curves}


Let $\cC$ be a balanced twisted curve over an algebraically closed field, $\cE$ a coherent sheaf on $\cC$, that is locally free at the nodes of $\cE$ (for the applications needed in this paper, the case of a locally free sheaf is sufficient: however, the added generality helps with the proof). Call $\pi \colon \cC \to C$ the coarse curve, $p_{1}$, \dots,~$p_{n}$ the marked points on $C$. For each $i$ call $r_{i}$ the index of $\cC$ at $p_{i}$. For each $i$ choose a section $p_{i} \to \Sigma_{i}$ of the marking gerbe $\Sigma_{i} \to p_{i}$, inducing an isomorphism $\Sigma_i \simeq \cB\bmu_{r_i}$. If $\cE$  is locally free at $p_{i}$, we have defined above the age of $\cE$ at $p_{i}$ as $   \age_{p_{i}}(\cE) = \age(\cE\mid_{\Sigma_{i}})$.

In the general case, when $\cE$ has torsion, this is not the correct definition. Consider the embedding
   \[
   \iota_{i} \colon \cB\bmu_{r_{i}} \cong \Sigma_{i} \into \cC;
   \]
it induces a pullback in the K-theory of coherent sheaves of finite projective dimension
   \[
   \iota_{i}^{*} \colon \mathrm{K}_{0}(\cC) \longrightarrow
      \mathrm{K}_{0}(\cB\bmu_{r_{i}}) = \mathrm{R}\bmu_{r_{i}}
   \]
via the usual formula
   \[
   \iota_i^{*}\cE = [\cE \otimes_{\cO_{\cC}}\cO_{\Sigma_{i}}]
      - [\tor_{1}^{\cO_{\cC}}(\cE, \cO_{\Sigma_{i}})];
   \]
the correct general definition is
   \[
   \age_{p_{i}}(\cE) = \age(\iota_{i}^{*}\cE).
   \]
This gives an additive homomorphism
   \[
   \age_{p_{i}} \colon \mathrm{K}_{0}(\cC) \longrightarrow \QQ.
   \]

By $\chi(\cE)$ we denote as usual the Euler characteristic of $\cE$ on $\cC$. For each $i$ we have $\H^{i}(\cC, \cE) = \H^{i}(C, \pi_{*}\cE)$, so the Euler characteristic is finite. In particular, since $\pi_{*}\cO_{\cC} = \cO_{C}$, we have
   \[
   \chi(\cO_{\cC}) = 1 - g
   \]
where $g$ is the arithmetic genus of $C$.

We define the \emph{degree} of $\cE$ as follows. First assume that $\cC$ is smooth and irreducible. Take an ordinary connected smooth curve $D$ with a finite morphism $\phi \colon D \to \cC$ (it is not hard to see that this exists). Call $d$ the degree of $\phi$. Then we set
   \[
   \deg_{\cC} \cE = \frac{1}{d}\deg_{D}\phi^{*}\cE.
   \]

It is easy to see that the degree is independent of the choice of $\phi$: if $\phi' \colon D' \to \cC$ is another choice, call $D''$ a component of the normalization of the fibered product $D \times_{\cC} D'$. This dominates both $D$ and $D'$, and then one uses standard properties of the degree.

If $\cC$ is not irreducible, take the normalization $\nu \colon
\overline{\cC} \to \cC$, pull back $\cE$, and sum the degrees over
the irreducible components of $\overline{\cC}$.

The degree is a rational number. If $\cE$ is locally free, then it
follows from the definition that the degree of $\cE$ is also the
degree of $\det \cE$, as usual; and we also have the formula
   \[
   \deg_{\cC} \cE = c_{1}(\cE) \cdot [\cC] = \int_{\cC}c_{1}(\cE) \cdot [\cC].
   \]

\begin{theorem} We have
   \[
   \chi(\cE) = \rk(\cE)\chi(\cO_{\cC}) + \deg \cE
      - \sum_{i=1}^{n} \age_{p_{i}}(\cE).
   \]
\end{theorem}

It is worth noticing that the stack structure at the nodes does not
intervene in the formula. This is due to the fact that the curve is
balanced.

This theorem can be deduced from To\"en's Riemann-Roch
theorem for stacks, but as it is not too much harder, we give a
direct argument.

\begin{proof}
In the following we will use the fact that Euler characteristic,
rank, degree and age are all additive invariants in
$\mathrm{K}_{0}(\cC)$.

Assume that $\cC$ is smooth.

First of all, assume that there exists a coherent sheaf $\cF$ on $C$
such that $\pi^{*}\cF = \cE$. Then the adjunction homomorphism $\cF
\to \pi_{*}\pi^{*}\cF$ is an isomorphism, because $\pi$ is flat.
Then $\chi(\cE) = \chi(\cF)$, $\deg_{\cC} \cE = \deg_{C}\cF$, and
$\age_{p_{i}}\cE = 0$ for all $i$; hence the formula follows from
ordinary Riemann--Roch applied to $\cF$.

The kernel and cokernel of the adjunction homomorphism
   \[
   \pi^{*}\pi_{*}\cE \longrightarrow \cE
   \]
are torsion, supported at the marked points of $\cC$; hence by additivity it is enough to prove the formula for torsion sheaves supported in the stack locus. Each such sheaf is an extension of sheaves of the form $\iota_{i*}L_{k}$, where $0 \leq k \leq r_{i}-1$, $\iota_{i}\colon \cB\bmu_{r_{i}}$ is the inclusion, and $L_{k}$ is the 1-dimensional representation of $\bmu_{r_{i}}$ defined by the character $\bmu_{r_{i}} \to \GG_{\mathrm{m}}$, $t \mapsto t^{k}$. So it suffices to prove the formula for the sheaves $\iota_{i*}L_{k}$.

It is easy to see that
   \[
   \deg \iota_{i*}L_{k} = \frac{1}{r_{i}}
   \]
for any $k$. Also we see that
   \[
   \chi(\iota_{i*}L_{k}) =
   \begin{cases}
   1 & \text{if } k = 0\\
   0 & \text{if } k \neq 0.
   \end{cases}
   \]
To complete the calculation let us compute $\age_{p_{i}} \iota_{i*}L_{k}$. Let $\cI_{\Sigma_{i}}$ be the sheaf of ideals of $\Sigma_{i}$ in $\cC$; by with tensoring $\iota_{i*}L_{k}$ the sequence
   \[
   0  \longrightarrow \cI_{\Sigma_{i}}  \longrightarrow \cO_{\cC}
      \longrightarrow \cO_{\Sigma_{i}}  \longrightarrow 0
   \]
we obtain that
   \[
   \iota_{i}^{*}\iota_{i*}L_{k} = [L_{k}]
      - [\mathrm{T}^{\vee}_{\Sigma_{i}} \otimes L_{k}]
   \]
where $\mathrm{T}^{\vee}_{\Sigma_{i}}$ is the cotangent space to $\Sigma_{i}$ in $\cC$ (this is a very particular case of the self-intersection formula in K-theory). But by definition of the isomorphism $\Sigma_{i} \cong \cB\bmu_{r_{i}}$, the tangent space to $\Sigma_{i}$ is $L_{1}$, hence $\mathrm{T}^{\vee}_{\Sigma_{i}}$ is $L_{r_{i}-1}$. From this we obtain the remarkable formula
   \[
   \age_{p_{i}} \iota_{i*}L_{k} =
   \begin{cases}\displaystyle
   -\frac{r_{i}-1}{r_{i}} & \text{if } k=0\\\displaystyle
   \frac{1}{r_{i}}        & \text{if } k\neq 0.
   \end{cases}
   \]
Now we plug in the values of the invariants and check the formula for the sheaves $\iota_{i*}L_{k}$. This completes the proof when $\cC$ is smooth.

In the general case, let $q_{1}$, \dots,~$q_{s}$ be the nodes of $C$. Call $l_{i}$ the index of $\cC$ at the node $q_{i}$; call $\Theta_{i}$ the residual gerbe of $\cC$ at the unique point living over $q_{i}$. Consider the normalization $\nu\colon \overline{\cC} \to \cC$; the moduli space $\overline{C}$ of $\overline{\cC}$ is the normalization of $C$. Call $q'_{i}$ and $q''_{i}$ the two inverse images of $q_{i}$ in $\overline{C}$, and call $\Theta'_{i}$ and $\Theta''_{i}$ the residual gerbe of  $\cC$ over the unique point over $q'_{i}$ and $q''_{i}$ respectively. The natural morphisms $\Theta_{i} \to \Theta'_{i}$ and $\Theta_{i} \to \Theta''_{i}$ are isomorphisms of gerbes. Furthermore, $\Theta'_{i}$ and $\Theta''_{i}$ are isomorphic to $\cB\bmu_{l_{i}}$, where the isomorphisms are chosen so that the action of $\bmu_{l_{i}}$ on the tangent space to $\cC$ at the gerbe is given by the embedding $\bmu_{l_{i}} \subseteq \GG_{\mathrm{m}}$. We give $\overline{\cC}$ the structure of a twisted curve by using these gerbes as a marking.

The fact that $\cC$ is balanced has the following important consequence: the two bandings by $\bmu_{l_{i}}$ on $\Theta_{i}$ that we obtain via the two isomorphisms $\Theta_{i} \cong \Theta'_{i}$ and $\Theta_{i} \cong \Theta''_{i}$ are opposite. Consider the pullback $\nu^{*}\cE$; its restrictions to $\Theta'_{i}$ and $\Theta''_{i}$ are isomorphic to the restriction of $\cE$ to $\Theta_{i}$, and they give dual representations of $\bmu_{l_{i}}$. This implies that
   \[
   \age_{q'_{i}} \nu^{*}\cE + \age_{q''_{i}} \nu^{*}\cE = c_{i},
   \]
where $c_{i}$ is the codimension of the space of invariants of the restriction of $\nu^{*}\cE$ to $\Theta'_{i}$ or $\Theta''_{i}$.

We have an exact sequence
   \[
   0 \longrightarrow \cO_{\cC} \longrightarrow \nu_{*}\cO_{\overline{\cC}}
     \longrightarrow \bigoplus_{i=1}^{s}\cO_{\Theta_{i}} \longrightarrow  0
   \]
from which we deduce that
   \[
   \chi(\cO_{\overline{\cC}}) = \chi(\nu_{*}\cO_{\overline{\cC}})
      = \chi(\cO_{\cC}) + s.
   \]
By tensoring with $\cE$, and keeping in mind that $\cE$ is locally free at the nodes of $\cC$, we get a sequence
   \[
   0 \longrightarrow \cE \longrightarrow \nu_{*}\nu^{*}\cE
     \longrightarrow \bigoplus_{i=1}^{s}(\cE \otimes_{\cO_{\cC}} \cO_{\Theta_{i}})
     \longrightarrow  0.
   \]
By taking invariants, and then using Riemann--Roch on the smooth twisted curve $\overline{\cC}$ together with the formulas above, we get
   \begin{align*}
   \chi(\cE) &= \chi(\nu_{*}\nu^{*}\cE)
      - \sum_{i=1}^{s} \dim_{k}\H^{0}(\cE \otimes_{\cO_{\cC}} \cO_{\Theta_{i}})\\
   &= \chi(\nu^{*}\cE)
      - \sum_{i=1}^{s}(\rk\cE - c_{i})\\
   &= \rk(\cE)\chi(\cO_{\overline{\cC}}) + \deg(\nu^{*}\cE)
      - \sum_{i=1}^{n}\age_{p_{i}}(\cE)
      - \sum_{i=1}^{s}(\age_{q'_{i}} \nu^{*}\cE + \age_{q''_{i}} \nu^{*}\cE)\\
      &\qquad- \sum_{i=1}^{s}(\rk\cE -
      \age_{q'_{i}} \nu^{*}\cE - \age_{q''_{i}} \nu^{*}\cE)\\
   &= \rk(\cE)\bigl(\chi(\cO_{\cC}) + s\bigr) + \deg(\cE)
      - \sum_{i=1}^{n}\age_{p_{i}}(\cE) - s \rk \cE\\
   &= \rk(\cE)\chi(\cO_{\cC}) + \deg(\cE) - \sum_{i=1}^{n}\age_{p_{i}}(\cE).
   \end{align*}

This concludes the proof.\end{proof}

\subsection{The \stringy Chow group and its grading}

Let $\cX$ be a smooth Deligne--Mumford stack as in the introduction. We define the rational \stringy Chow group
 of $\cX$ to be  $$\stChow{\cX}: = A^*(\riX)_\QQ.$$ We make this into a graded group using the following rule (see \cite{Chen-Ruan}, \cite{DHVW}, \cite{Zaslow}):
 $$\stChow[a]{\cX} = \oplus_\Omega A^{a-\age(\Omega)}(\Omega)_\QQ, $$
 where the sum is taken over all connected components $\Omega$ of
 $\riX$.

\subsection{The small quantum Chow ring}

The small quantum Chow ring is an algebra over the completed monoid-algebra  $\QQ\llbracket  N^+(X) \rrbracket$, where we denote the monomial corresponding to a class $\beta$ by $q^\beta$. As a group, the quantum Chow ring is $$\QA^* (\cX) : = \stChow{\cX} \llbracket N^+(X)\rrbracket.$$
We define a product  on $\QA^* (\cX)$ by specifying the product of monomials, as follows:

$$ \gamma_1 * \gamma_2 = \sum_{\beta\in N^+(X)} \left\langle \gamma_1,\gamma_2,*\right\rangle_{0,\beta}\, q^\beta.$$

We define the \emph{degree} of $q^{\beta}$ to be
$\beta \cdot c_{1}(\cT_{\cX})$. 

\begin{theorem}
The product defined above makes $\QA^* (\cX)$ into a commutative, associative ``pro-$\QQ$-graded'' ring, in the sense that the product of two homogeneous elements of degrees $a,b$ is homogeneous of degree $a+b$.
\end{theorem}

\begin{proof}
Commutativity is immediate. Associativity is, as usual, a consequence of Theorem \ref{Thm:WDVV}. It remains to check the claim about grading.

Consider the summand $ \left\langle \gamma_1,\gamma_2,*\right\rangle_{0,\beta}\, q^\beta$ in the formula above, where we assume that $\gamma_{1}$ and $\gamma_{2}$ are each supported on a single component $\Omega_{1},$ respectively $\Omega_{2}\subset\riX$, of corresponding ages $a_{1}, a_{2}$. We need to show that it has degree $$\deg \gamma_{1}+\deg \gamma_{2} = (\codim_{\Omega_{1}} \gamma_{1} +{a_{1}}) +  (\codim_{\Omega_{2}} \gamma_{2} +{a_{2}}).$$ Similarly, it is enough to calculate the degree of a single term
   \[
   \left\langle \gamma_1,\gamma_2,*\right\rangle_{0,\beta,\Omega_{3}}\, q^\beta
   \]
of $ \left\langle \gamma_1,\gamma_2,*\right\rangle_{0,\beta}\, q^\beta$ lying in the Chow group of a component $\Omega_{3}\subset \riX$ having age $\check a_{3}$.

Given a stable map $f: \cC \to \cX$ corresponding to a geometric point of a component $\cK\subset \cK_{0,3}(\cX,\beta)$ with evaluations $e_{1}:\cK \to \Omega_{1}$, $e_{2}:\cK \to \Omega_{2}$ and $\check e_{3}:\cK \to \Omega_{3}$, the bundle $f^{*}\cT_{\cX}$ has ages $a_1, a_2$ and $a_{3}$ at the three markings, with $$a_{3} + \check a_{3} = \dim \cX - \dim \Omega_{3}.$$ This is explained in \cite{Chen-Ruan}, the point being that $\dim_{\xi,\alpha} \Omega_{3} = \mathrm{rank} (\cT_{\cX}^{\mu_{r(\Omega_{3})}}) = \dim \cX - (a_{3}+ \check a_{3})$, as noted in the last section.  We denote the class of $f_{*}\cC$ by $\beta'$, and clearly we have $\pi_{*}\beta' = \beta$.

We can now calculate dimensions at $f:\cC \to \cX$ using Riemann-Roch: first, the dimension of $[\cK]^{\vir}$ is given by
$$ \chi (f^{*}\cT_{\cX}) = c_{1}(\cT_{\cX})\cdot \beta' + \dim \cX - a_{1}-a_{2}-a_{3}.$$
Therefore we have
\begin{align*}
   \dim((e_{1}^{*} \gamma_{1}\cup e_{2}^{*} \gamma_{2})\cap [\cK]^{\vir})
 &= \deg q^{\beta} + \dim \cX -a_{1}-a_{2}-a_{3}\\
 &\phantom{= \deg q^{\beta} + \dim \cX}- \codim_{\Omega_{1}} \gamma_{1} - \codim_{\Omega_{2}} \gamma_{2} \\
 &= \deg q^{\beta} + \dim \cX -\deg \gamma_{1}- \deg \gamma_{2} - a_{3}.
\end{align*}
Pushing forward by $\check e_{3}$ we obtain
\begin{align*}
   \codim &\left(\check e_{3\,*}    ((e_{1}^{*} \gamma_{1}\cup e_{2}^{*}) \gamma_{2})\cap [\cK]^{\vir})\right) \\
    &= \dim\Omega_{3} - (\deg q^{\beta} + \dim \cX -\deg \gamma_{1}- \deg \gamma_{2} - a_{3} ) \\
   &= \deg \gamma_{1}+ \deg \gamma_{2} + (a_{3} +\dim \Omega_{3} - \dim \cX) - \deg q^{\beta} \\
   &= \deg \gamma_{1}+ \deg \gamma_{2} - \check a_{3}- \deg q^{\beta}
\end{align*}
 and therefore
    \begin{align*}
   \deg \left(\check e_{3\,*}    ((e_{1}^{*} \gamma_{1}\cup e_{2}^{*}) \gamma_{2})\cap [\cK]^{\vir})\right) = \deg \gamma_{1}+ \deg \gamma_{2} - \deg q^{\beta}.
   \end{align*}
   It follows that
  \begin{align*}
   \deg \left(\check e_{3\,*}    ((e_{1}^{*} \gamma_{1}\cup e_{2}^{*}) \gamma_{2})\cap [\cK]^{\vir})q^{\beta}\right)  = \deg \gamma_{1}+ \deg \gamma_{2},
   \end{align*}
as required.
\end{proof}

We note that, while the grading has rational degrees, the
denominators which appear are bounded. The fact that the $q^{\beta}$
have degrees with bounded denominators follows from the fact that the
ring is finitely generated, and more explicitly through
Proposition~\ref{Prop:bound-denominators}; and the denominators in
the ages of a connected component are bounded by the exponent of the
automorphism group of a geometric point of $\cX$.

An identical construction using the topological Gromov-Witten classes, gives us a definition
of the small quantum cohomology ring $\QH^*(\cX)$. This is an interesting ring
structure on $\H^*(\riX)\llbracket N(X)^+\rrbracket$. Its associativity follows from Theorem \ref{Thm:WDVV-homology}.

\subsection{Stringy Chow ring}

By setting the $q$'s to zero in $\QA^*(\cX)$, we get an interesting product structure on $\stChow{\cX}$ which we
refer to as the stringy Chow ring.  This is the Chow analogue of the
\emph{orbifold cohomology} or \emph{Chen-Ruan cohomology ring}, which arises by doing exactly the same thing to the small quantum cohomology ring of $\cX$.  In \cite{AGV} we show that it is possible to define these products with integral coefficients, provided one works with $A^*(\iX)$ instead of $A^*(\riX)$.

\subsection{The big quantum cohomology ring}

As in the Gromov-Witten theory of a manifold, the full topological WDVV equation,
Theorem \ref{bigwdvv}, is
essentially equivalent to the associativity of a big quantum cohomology ring
whose multiplication is defined in terms of the basis $\{
\alpha_i\}$ of $\H^*(\riX)$ as
$$\mu * \nu\  \ =\ \  \sum_{ij} \sum_n \ \frac{1}{n!}\ \big\langle \mu, \nu, T^n, \alpha_i \big\rangle_\beta\  q^\beta\
\tilde g^{ij}\ \alpha_j .$$ Here $T=\sum \alpha_i x_i$ should be
expanded formally, so that the resulting product is a power series
in the $x_i$ and $q_i$ whose coefficients record every genus zero
Gromov-Witten invariant including the classes $\mu$ and $\nu$.  The same arguments as in the previous subsection show
that this gives rise to a pro-$\QQ$-graded, associative, commutative ring structure on
$\H^*(\riX)\llbracket N^+(X)\rrbracket\llbracket x_1,\ldots , x_M\rrbracket$.

\section{A few useful facts}

We collect here some useful facts analogous to standard properties of Gromov-Witten
invariants for manifolds.  While this will not constitute an exhaustive list, we hope that readers will see how standard facts and techniques from Gromov-Witten theory carry over to this context.  Most sources of difference come from the need to systematically replace $\cX$ with $\riX$ at various points in formulating the theory, or from a difference in the relationship between the spaces of twisted stable maps and their universal curves, which we explain now.

\subsection{Universal curve}

A critical fact in ordinary Gromov-Witten theory is that $\ocM_{g,n+1}(X,\beta)$ is the
universal curve over $\ocM_{g,n}(X,\beta)$.  While this is not true for the orbifold theory,
we have a similar result.

\begin{proposition}\label{universalcurve} The
universal curve over $\cK_{g,n}(\cX,\beta)$ is naturally identified
with the open and closed substack $\cU \subseteq \cK_{g,n+1}(\cX,\beta)$ for
which the $(n+1)$st marked point is untwisted.  Moreover, if we
consider the flat morphism $\pi: \cU \to \cK_{g,n}(\cX,\beta)$, then
we have an equality of virtual fundamental classes $[\cU]^\vir =
\pi^*[\cK_{g,n}(\cX,\beta)]^\vir$, where $[\cU]^\vir$ denotes the
restriction of $[\cK_{g,n+1}(\cX,\beta)]^\vir$ to $\cU$.
\end{proposition}
Another way to say this, which is useful in practice, is that the universal curve $\cU$
fits in the cartesian square:
$$\xymatrix{\cU\ar[r]\ar[d]& \cX\ar[d]^i\\
\cK_{g,n+1}(\cX,\beta)\ar[r] & \riX}
$$
where $i$ denotes the inclusion of $\cX\cong \riX[1]$ into the inertia stack.

\begin{proof}
The identification of the universal curve is essentially Corollary 9.1.3 in \cite{AV}.  The identification of the virtual classes follows by the same reasoning as the analogous statement in \cite{Behrend}.
\end{proof}

\subsection{Degree zero invariants}
The identification of the $(n+1)$-pointed space with the universal curve in ordinary Gromov-Witten theory implies that the degree zero invariants carry very little information.  In particular, it implies that degree zero, genus zero invariants all vanish with the exception of the three point invariants, which simply compute the classical cohomology ring.  This is not at all true for the orbifold theory, which is precisely why the ring structure on orbifold cohomology is interesting.
However, one consequence of this fact does continue to hold.  If we let $1$ denote the fundamental class of the identity component of $\riX$, then we have the following fact.

\begin{proposition} The Gromov-Witten number $\langle 1, \delta_1, \ldots , \delta_n \rangle_{g,\beta}$
vanishes unless $g=0$ and $n=2$, in which case we have
$$\langle 1, \delta_1, \delta_2 \rangle_{0,0}= \int_{\riX} \frac{1}{r} \delta_1 \cup \iota^*\delta_2.$$

\end{proposition}

\begin{proof}
Given Proposition \ref{universalcurve}, the vanishing portion follows from the same arguments as in ordinary Gromov-Witten theory (see \cite{KM}).  The precise formula comes from the easy fact that $e_1^{-1}(\riX[1])$ in $\cK_{0,3}(\cX,0)$ is naturally isomorphic to $\iX$ and we can identify $e_2$ with the standard rigidification map.  Using this identification, $e_3$ is then identified with the composition of the rigidification map with $\iota$ and the result follows.
\end{proof}

This has the immediate consequence that $1$ is the identity element in both the small and big quantum cohomology rings.  (A completely analogous lemma shows that 1 is the identity in the quantum Chow ring.)

\subsection{Gravitational descendants}
In studying higher genus Gromov\ddash Witten theory it is important to include the descendant classes.  These are analogues
of the Mumford--Morita--Miller classes on $\overline{M}_{g,n}$ involving the Chern
classes of the normal bundles of  the $n$ sections.

One can define these  classes in the orbifold theory in much the
same way as in the ordinary theory.  On $\cK_{g,n}(\cX,\beta)$, one defines
 $n$ tautological line bundles, $\cL_i$.  There is more than one
way to define these, but the most straightforward is to take $\cL_i$ to
be the bundle whose fiber at a point is the cotangent space to the
corresponding \emph{coarse} curve at the $i$th marked point, in other words, the pullback via the $i$-th section on the universal \emph{coarse} curve of the sheaf of relative differentials.  We
remark that while one might like to use the cotangent space to the
twisted curve, these will not form a line bundle on
$\cK_{g,n}(\cX,\beta)$, but only on the gerbe corresponding to the
$i$th marking.  However, one could decide to push down the Chern
class of that ``twisted" bundle, which would simply be $1/r$ times the
Chern class of the $\cL_i$ as we are defining them here.  We will not
use that convention.

Let $\psi_i = c_1(\cL_i)$.   Then given classes $\delta_1, \ldots
\delta_n$ in $\H^*(\riX)$, we define the following invariants:
$$\langle \, \tau_{a_1}(\delta_1) \cdots
\tau_{a_n}(\delta_n)\, \rangle_{g,\beta} =
\int_{\cK_{g,n}(\cX,\beta)]^\vir} e_1^*(\delta_1) \cup \psi_1^{a_1}
  \cup \cdots \cup e_n^*(\delta_n)\cup \psi_n^{a_n}.$$

For simplicity of statements we  define any class involving $\tau_{-1}$ in the formulas below to be 0.

With these definitions, the standard equations among descendants for manifolds (cf. \cite[Section 1.2]{Pandharipande}) hold in the orbifold setting with no changes (keeping in mind that 1 means the fundamental class of the identity component of the inertia stack).

\begin{theorem}
Assume $(\beta,g,n)$ is not any of $\beta=0, g=0,
n<3$ or $\beta=0,g=1,n=0$. Then
\begin{enumerate}
\item
(Puncture or String Equation)
\begin{align*}
\langle\, \tau_{a_1}(\gamma_1)& \cdots \tau_{a_n}(\gamma_n)\tau_0(1)\,\rangle_{g,\beta} \\&=
\sum_{i=1}^n \langle\, \tau_1(\gamma_1) \cdots \tau_{a_{i-1}}(\gamma_
{i-1})\tau_{a_i-1}(\gamma_i)
\tau_{a_{i+1}}(\gamma_{i+1}) \cdots \tau_{a_n}(\gamma_n)\,\rangle_{g,\beta}
\end{align*}
\medskip

\item (Dilaton Equation)
$$
\langle\, \tau_{a_1}(\gamma_1) \cdots \tau_{a_n}(\gamma_n)\tau_1(1)\,\rangle_{g,\beta} =
(2g-2+n)\langle\,\tau_{a_1}(\gamma_1) \cdots \tau_{a_n}(\gamma_n)\,\rangle_{g,\beta}
$$
\medskip

\item (Divisor Equation)  For $\gamma$ in $\H^2(X) \subset \H^2_{\mathrm{orb}}(X)$ (but not for an
arbitrary element of $\H^2_{\mathrm{orb}}(X)$), we have
\begin{align*}
\langle\, \tau_{a_1}(\gamma_1)&\cdots
\tau_{a_n}(\gamma_n)\tau_0(\gamma)\,{\rangle}_{g,\beta}
\\=&
\left(\int_\beta \gamma\right) \cdot \langle \, \tau_{a_1}(\gamma_1) \cdots \tau_
{a_n}(\gamma_n)\,\rangle_{g,\beta}
\\&+ \sum_{i=1}^n \langle\, \tau_1(\gamma_1) \cdots
\tau_{a_{i-1}}(\gamma_{i-1})\tau_{a_i-1}(\gamma_i\cup \gamma)
\tau_{a_{i+1}}(\gamma_{i+1}) \cdots \tau_{a_n}(\gamma_n)\,\rangle_{g,\beta}
\end{align*}
\end{enumerate}
\end{theorem}


\begin{proof}
We reduce to the untwisted case. In all these equations
the intersection happens on the open and closed substack $\cU \subset
\cK_{g,n+1}(\cX, \beta)$ of Proposition  \ref{universalcurve}. The key
commutative diagram is the following:
$$\xymatrix{
\cU \ar[r]^{\mu_{n+1}} \ar[d]_{\pi_\cK}
                                 & \ocM_{g,n+1}(X, \beta)\ar[d]^{\pi_\cM} \\
 \cK_{g,n}(\cX, \beta)\ar[r]^{\mu_{n}}
                                 & \ocM_{g,n}(X, \beta)
}$$
This is not quite a fiber diagram, but $\cU$ has the same moduli space
as the fibered product, so the projection formula still holds. If
$\alpha$ is a cohomology class on $\cK_{g,n}(\cX, \beta)$ and $\psi$
is a cohomology class on $\ocM_{g,n+1}(X, \beta)$, we have by
proposition \ref{universalcurve}
$$ \int_{[\cU]^\vir} \pi_\cK^*\alpha\cup\mu_{n+1}^*\psi \ \ = \ \
\int_{[\cK_{g,n}(\cX, \beta)]^\vir} \alpha \cup
\mu_n^*\pi_{\cM\,*}\psi.$$ So any identity which holds for
$\pi_{\cM\,*}\psi$ (keeping in mind that $X$ may be singular) can be
used here.

For all the equations use $\alpha = \prod_{i=1}^n e_i^*\gamma_i$. For
the Puncture Equation use $\psi = \prod_{i=1}^{n} \psi_i^{a_i}$ and
the equation
$\pi_{\cM\,*}\psi = \sum_{i=1}^n(\prod_{i=1}^{n}
\psi_i^{a_i})/\psi_i$.  For the Dilaton Equation use  $\psi =
(\prod_{i=1}^{n} \psi_i^{a_i})\psi_{n+1}$  with $\pi_{\cM\,*}\psi =
(2g-2+n) \prod_{i=1}^{n} \psi_i^{a_i}$.  For  the  Divisor equation
use
$\psi =   (\prod_{i=1}^{n} \psi_i^{a_i}) \cup
e_{\cM_{n+1},n+1}^*\gamma$, with the equation $\pi_{\cM\,*}\psi =
\int_\beta \gamma \cdot \prod_{i=1}^{n} \psi_i^{a_i} + \sum_{i=1}^n
(\prod_{i=1}^{n}\psi_i^{a_i})/\psi_i \cap e_{\cM_{n},i}^*\gamma,$
noticing that  $e_i^*\gamma_i \cup u_n^* e_{\cM, i}^* \gamma =
e_i^*(\gamma_i\cup \gamma)$.
\end{proof}

\begin{remark} There is also a Topological Recursion Relation valid in this context, treated in \cite{Tseng}, section 2.5.5.
\end{remark}

\section{An example: the weighted projective line}\label{Sec:example}

We conclude by giving a nontrivial
calculation of the quantum Chow
ring of a stack.
We consider one of the simplest possible classes of examples -- the
weighted projective lines.
We fix two positive integers, $a$ and $b$, and consider the one
dimensional weighted projective space $\cX=\proj(a,b)$ which is the
stack quotient of a punctured two dimensional affine space by the
action of $\GG_m$ with weights $a$ and $b$.  The coarse moduli space
of $\proj(a,b)$ is always $\proj^1$.  If $a$ and $b$ are relatively
prime, then $\proj(a,b)$ is a twisted curve.  Otherwise this stack
has a generic stabilizer.  We will denote by $0$, the point with
stabilizer $\mu_a$ and by $\infty$ the point with stabilizer
$\mu_b$.

For the convenience of the reader, we collect the basic facts about
this stack here.  A morphism from a scheme $Z$ to $\proj(a,b)$ is
given by choosing a line bundle $L$ on $Z$ together with sections
$s_1 \in \Gamma(Z,L^{\otimes a})$ and $s_2 \in \Gamma(Z,L^{\otimes
b})$ with no common zeroes.  2-morphisms are given by morphisms
between line bundles which take the sections to the sections.  Note
that if $a=b=1$ we get the usual description of $\proj^1$ and there
are no nontrivial 2-automorphisms.  By descent, we get the same
description of maps from a stack to $\proj(a,b)$.  We let $\oh(1)$ denote the
line bundle on $\proj(a,b)$ corresponding to the identity
morphism.  Then $\Pic
(\proj(a,b)) = \ZZ \oh(1)$ and there are sections of $\oh(a)$ and
$\oh(b)$ vanishing at $\infty$ and 0 respectively.  The degree of
$\oh(1)$ is $\frac{1}{ab}$.  Finally, $T\proj(a,b) \cong \oh(a+b)$.

The inertia stack of this stack is straightforward to describe. Note
that since $\GG_m$ is an abelian group, the quotient presentation of
$\cX$ endows the inertia group of each point of $\proj(a,b)$ with an
embedding in $\GG_m$.  (This character of the isotropy group is also
its action on the fiber of $\oh(1)$.)  Because of this, each
irreducible component of $I_{\mu_r}\cX$ is canonically associated with
the unit in $\ZZ/r\ZZ$ which relates that fixed embedding to the one
obtained by composition with the embedding of $\mu_r$ in $\GG_m$.
Below we will use the obvious convention which identifies the set of
elements of $\ZZ/d\ZZ$ with the set of units in $\ZZ/k\ZZ$ for all
positive integers $k$ dividing $d$.

Let $d=gcd(a,b)$. For each element of $\ZZ/d$, there is a one
dimensional component of $I_\mu\cX$ which is isomorphic to
$\proj(a,b)$. For each element of $\ZZ/a\ZZ$ which is not divisible
by $a/d$, there is a zero dimensional component of $I_\mu(\cX)$
lying over the point $0$. Each of these components is isomorphic to
$\cB\mu_a$. Similarly, for each element of $\ZZ/b\ZZ$ not divisible by
$b/d$ there is a component lying over the point $\infty$ which is
isomorphic to $\cB\mu_b$.  We hope that the following picture of the
inertia stack of $\proj(4,6) \cong \ocM_{1,1}$ will make this labeling system clear:



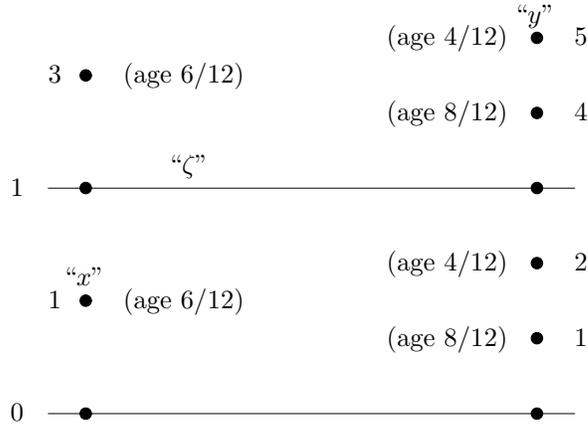
\begin{center}
\setlength{\unitlength}{1mm}
\begin{picture}(80,55)(0,0)
\put( 5,14){1}
\put(10,15){\blacken\ellipse{1.5}{1.5}}
\put(15,14){(age 6/12)}
\put( 7,17){``$x$''}
\put( 5,44){3}
\put(10,45){\blacken\ellipse{1.5}{1.5}}
\put(15,44){(age 6/12)}
\put(50, 9){(age 8/12)}
\put(70,10){\blacken\ellipse{1.5}{1.5}}
\put(75, 9){1}
\put(50,19){(age 4/12)}
\put(70,20){\blacken\ellipse{1.5}{1.5}}
\put(75,19){2}
\put(50,39){(age 8/12)}
\put(70,40){\blacken\ellipse{1.5}{1.5}}
\put(75,39){4}
\put(50,49){(age 4/12)}
\put(70,50){\blacken\ellipse{1.5}{1.5}}
\put(75,49){5}
\put(67,52){``$y$''}

\path(5,0)(75,0)
\put(10,0){\blacken\ellipse{1.5}{1.5}}
\put(70,0){\blacken\ellipse{1.5}{1.5}}
\put(0,-1){0}
\path(5,30)(75,30)
\put(10,30){\blacken\ellipse{1.5}{1.5}}
\put(70,30){\blacken\ellipse{1.5}{1.5}}
\put(0,29){1}
\put(21,32){``$\zeta$''}
\end{picture}
\end{center}

To give a presentation of the quantum Chow ring we need to choose
generators.  A convenient way to
make this choice is as follows. Choose integers $m$ and $n$ such
that $ma+nb=d$.  Set $A=a/d$ and $B=b/d$. We take $\zeta$ to be the
fundamental class of the one dimensional component of $\riX$ corresponding to
$1\in \ZZ/d\ZZ$ (if $d=1$, take $\zeta=1$), we let $x$ be the fundamental class of the
component lying over $0$ which corresponds to $n\in \ZZ/a\ZZ$, and
we let $y$ be the fundamental class of the component lying over
$\infty$ corresponding to $m\in \ZZ/b\ZZ$. In the example of
$\PP(4,6)$ above, we chose $n=1, m=-1$ and indicated the components where the resulting
$x,y$ as well as $\zeta$ serve as fundamental classes. (If $d=a$ or $d=b$ some
of these zero dimensional components don't exist.  We will ignore
this case in what follows, but the results all hold with essentially
identical proofs if we take $x$ or $y$ to be the fundamental class
of the appropriate zero dimensional substack of the inertia stack
associated with $m$ or $n$ in $\ZZ/d\ZZ$.)  One consequence of
choosing $x$ and $y$ in this manner is that they have minimal age.
Under this convention, we find that $\deg(x)=1/A$ and $\deg(y)=1/B$,
while $\deg(\zeta)=0$.
 Following reasoning similar to that at the end of \cite{AGV} it is
easy to calculate that the stringy Chow ring of $\proj(a,b)$ is
$\QQ[\zeta,x,y]\, /\, \langle xy, Ax^A-By^B\zeta^{n-m}, \zeta^d-1\rangle$. (Note that the factor  $\zeta^{n-m}$ in the second relation is missing in \cite{AV-families}.)

The N\'eron-Severi group of $\proj(a,b)$ has rank one, so to compute
the quantum Chow ring, we need to introduce one further generator,
$q$.  We normalize this by selecting our generator of
$N(\proj(a,b))$ to correspond to the minimal positive degree map
from a twisted curve to $\proj(a,b)$.  If $gcd(a,b)=1$, then this
minimal degree map can be taken to be the identity map, if not, then
the presence of a generic stabilizer of order $d$ forces any map
from a twisted curve to have degree divisible by $d$, we will see
that there does exist a map of degree exactly $d$.  In other words,
we take the generator of $N(\proj(a,b))$ to be $d$ times the
fundamental class.

Since we know that $\QH^*(\cX)$ is free as a $\QQ\llbracket q\rrbracket$ module, it
follows that there is a presentation for $\QH^*(\cX)$ of the form
$$\QQ\llbracket q\rrbracket[\zeta,x,y]\ \big/\ \langle R_1, R_2, R_3 \rangle$$
where we have
$$R_1 \equiv xy \mod q$$
$$R_2 \equiv Ax^A-By^B\zeta^{n-m} \mod q$$
$$R_3 \equiv \zeta^d-1 \mod q.$$

Since the degree of $q$ is $1/A + 1/B$, no monomial in the
generators containing a $q$ can possibly have degree equal to the
degrees of $R_2$ or $R_3$.  Hence all that remains to compute is the
quantum product of $x$ with $y$.  By degree considerations, the only
possibility for the form of $R_1$ is then
$$xy=q(c_0+c_1\zeta + \cdots +c_{d-1}\zeta^{d-1})$$ where the $c_i$
are rational numbers. We will establish that $c_0=1$ and that the
other $c_i$ all vanish. As $c_0$ is the coefficient of the
fundamental class in $x*y$, it is determined by considering the
moduli space of maps from a 3 pointed stacky $\proj^1$ which
actually has only two stacky points of indices $a$ and $b$.  We let
$C_{a,b}$ denote this curve. $\Pic(C_{a,b})$ is generated by line
bundles $L_0$ and $L_\infty$ of degrees $\frac 1a$ and $\frac 1b$
satisfying the single relation $L_0^{\otimes a} \cong
L_\infty^{\otimes b}$.  In particular, when $a$ and $b$ are not
relatively prime, so that $C_{a,b} \not\cong \proj(a,b)$ we find
that there is torsion in the Picard group. Also, the restriction map
$r_0: \Pic (C_{a,b}) \to \Pic(\cB\mu_a)$ takes $L_0$ to the standard
generator.  (A possibly confusing point here is that even when $a$
and $b$ are relatively prime so that $C_{a,b}\cong \proj(a,b)$, our
 identification of the isotropy groups of the substacks
supported at 0 with $\mu_a$ are different, since in $\proj(a,b)$ we
are using the action of $\mu_a$ on the fiber of $\oh(1)$ as the
standard representation, whereas on $C_{a,b}$ we are using the
action on the tangent space.)

We are looking for a map of degree $d$ from $C_{a,b}$ to
$\proj(a,b)$. Since we can compute the degree of a morphism $f:
C_{a,b} \to \proj(a,b)$ by comparing the degree of $\oh(1)$ to the
degree of $f^*\oh(1)$ we see that we need to find a line bundle $L$
of degree $\frac{d}{ab}$ such that $L^{\otimes a}$ has a section
vanishing at $\infty$ and $L^{\otimes b}$ has a section vanishing at
$0$.  If we denote this line bundle $L= L_0^{\otimes z_0} \otimes
L_\infty^{\otimes z_\infty}$ then in order for the first point to
evaluate to the correct component of the inertia stack, we need that
$z_0 \equiv n \mod a$ and for the second point we have the analogous
condition that $z_\infty \equiv m \mod b$.  It follows that we must
have $L\cong L_0^m\otimes L_\infty^n$.  Thus the relevant space of
morphisms from $C_{a,b}$ to $\proj(a,b)$ can be identified with the
space of pairs $s_1 \in \Gamma(L^{\otimes a}), s_2 \in
\Gamma(L^{\otimes b})$ modulo the action of $\CC^*$ acting by scalar
multiplication on $L$.  Since we are assuming
that $d\neq a$ and $d \neq b$, there is a unique section of both
$L^{\otimes a}$ and $L^{\otimes b}$ by degree considerations. We conclude that the space of 3 pointed maps with irreducible source curve is the quotient of $\GG_m \times \GG_m$ by the linear action of $\GG_m$ with
weights $a$ and $b$.  This is simply  $\GG_m \times \mu_d$.
There are two additional points where the source is reducible with a component collapsed over either zero or infinity.  The reader may verify that the full moduli space is isomorphic to $C_{a,b} \times \cB\mu_d$, but to compute the relevant pushforward, the exact structure of the compactification is irrelevant.

Since this space has the expected
dimension, the virtual fundamental class is the ordinary fundamental
class, which pushes forward to the fundamental class of $\proj(a,b)$
and we conclude that $c_0=1$.  To see that the other $c_i$
vanish we just need to verify that there are no representable
morphisms of minimal degree from a three pointed twisted genus zero
curve where the first two points are as before, but the third
twisted point has nontrivial stacky structure. To find such a map,
we would need to find an integer $D|d$, and a line bundle of degree
$\frac{d}{ab}$ on $C_{a,b,D}$ satisfying all of the conditions as
before as well as being nontrivial when restricted to $\cB\mu_D$. This
is obviously impossible.

We conclude that
$$\QA^*(\proj(a,b))\ \ =\ \  \QQ\llbracket q\rrbracket[\zeta,x,y]\ \big/\ \langle
xy-q,\, Ax^A-By^B\zeta^{n-m},\, \zeta^d-1\rangle.$$

\appendix
\section{Gluing of algebraic stacks along closed substacks} \label{Sec:Gluing}
\subsection{} We introduce a gluing construction for
Artin stacks.



\begin{proposition}\label{Prop:gluing}
Let $Z,X_1,X_2$ be algebraic stacks, and  assume given
 $$ \xymatrix{ *+<12pt,12pt>{Z} \ar@{^{(}->}[r]^{i_1} \ar@{^{(}->}[d]_{i_2} & X_1 \\
X_2  } $$
  where $i_1,i_2$ are closed embeddings. Then there exists
an algebraic stack $X$ and a diagram
 $$ \xymatrix{ *+<12pt,12pt>{Z} \ar@{^{(}->}[r]^{i_1} \ar@{^{(}->}[d]_{i_2} & X_1\ar[d] \\
X_2 \ar[r] & X} $$
such that the diagram is co-cartesian, namely, for any algebraic
stack $\cM$, the natural functor
$$\Hom(X,\cM) \to \Hom(X_1,\cM)\ \ \mathop{\times}\limits_{\Hom(Z,\cM)}\
\  \Hom(X_2,\cM)$$
is an equivalence of categories. Such $X$ is unique up to a unique isomorphism.
\end{proposition}

\begin{corollary}\label{Cor:limits}
Let $Z, Y$ be algebraic stacks and $i_1, i_2: Z \hookrightarrow Y$ closed embeddings with disjoint images.    Then there exists
\begin{enumerate}
\item an algebraic stack $X$,
\item a morphism $\pi: Y \to X$, and
\item a 2-isomorphism $\alpha: \pi\circ i_1 \to\pi\circ i_2$
\end{enumerate}
 such that $(X,\pi,\alpha) =
 \displaystyle\varinjlim( Z \double Y)$.
\end{corollary}

In other words, we can glue together the two copies of $Z$ in $Y$, obtaining $X$

\subsection{Gluing of schemes and algebraic spaces}\label{Sec:gluing}

 Given a scheme $Z$ together with a pair of closed embeddings
of $Z$ in schemes $X_1$ and $X_2$, we can define a new scheme $X_1\cup_Z X_2$.
It is determined by the universal property that a morphism
from $X_1 \cup_Z X_2$ to a scheme $W$ is given by a morphism from $X_1$
to $W$ and a morphism from $X_2$ to $W$ whose restrictions to $Z$
agree.  To construct this scheme, we can just do the construction
for affines, where this amounts to taking a fibered product of rings.

In the following we use the fact that a similar construction exists in
the category of algebraic spaces.
One way to prove
that it exists is to use our proof below, replacing ``algebraic
space'' by ``scheme'' and ``algebraic
stack'' by ``algebraic space''.

Our first lemma shows that the universal property of the gluing of
algebraic spaces is preserved in the 2-category of stacks.

\begin{lemma}\label{Lem:gluing-is-coproduct}  Let $\cM$ be an
algebraic stack, and let
$$ \xymatrix{ *+<12pt,12pt>{Z} \ar@{^{(}->}[r]^{i_1} \ar@{^{(}->}[d]_{i_2} & X_1\ar[d] \\
X_2 \ar[r] & X} $$
be a co-cartesian diagram of algebraic spaces, where $i_1,i_2$ are
closed embeddings.
Then the natural functor
$$\cM(X) \to \cM(X_1) \times_{\cM(Z)} \cM(X_2)$$
is an equivalence, where the fibered product is taken in the sense of categories.
\end{lemma}

\emph{Proof.} We construct a functor inverse to the given one. Let $R
\double U $ be  a presentation of $\cM$. Assume given an object of the fibered product on the right hand side, namely
\begin{enumerate}
\item objects $f_i\in \cM(X_i) $, and
\item an isomorphism $\alpha: i_1^* f_1 \to i_2^* f_2$.
\end{enumerate}
   Explicitly in terms of the presentation, we are given $U_i =
   f_i^*U$ and $R_i = f_i^*R$ and morphisms of groupoids
   $$ \xymatrix{ R_i \doubledown \ar[r] &R \doubledown \\
                 U_i \ar[r]             & U.         }
      $$
Moreover, $\alpha$ gives an isomorphism
$$i_1^* ( R_1 \double U_1) \to i_2^* ( R_2 \double U_2).$$
Since $R_i$ and $U_i$ are algebraic spaces, we can form a groupoid in
algebraic spaces $R_X \double U_X$ presenting $X$, by gluing $U_1,U_2$
along the
morphism $\alpha:i^*U_1 \to i^*U_2$, and similarly for $R_X$. We
obtain a morphism of groupoids
$$ \xymatrix{ R_X \doubledown \ar[r] &R \doubledown \\
                 U_X \ar[r]             & U         }
      $$
giving a morphism $X \to \cM$.
The construction of the functor on the level of arrows, and the fact
that the functors are inverses, is left to the reader. \qed

\subsection{Extension of atlases}
\begin{lemma} Let $Z\subset X$ be a closed embedding of
schemes. Assume $U_Z \to Z$ is a smooth morphism. Then there exists a
Zariski open covering $U_Z' \to U_Z$ and a smooth morphism $U_X' \to
X$ with  $U_X' \times_X Z \cong U_Z'$.
\end{lemma}

\begin{proof}
We may assume $U_Z,Z$ and $X$ are affine, and embed $U_Z
\subset \bbA^n_Z$ for some $n$.  The subscheme $U_Z$ is locally a complete
intersection in $\bbA^n_Z$, so for every point $u\in U_Z$ we can find
elements $(g_1,\ldots g_N)\subset \cO_Z[x_1,\ldots,x_n]$ and   a
Zariski neighborhood $V \subset \bbA^n_Z$ such that $U_Z^u :=
U_Z \cap V
\subset V$ is the complete intersection of the zero schemes of
$g_i$. Choose $\tilde g_i\subset \cO_X[x_1,\ldots,x_n]$ lifting
$g_i$, and let $W$ be the zero scheme of $(\tilde g_1,\ldots,\tilde
g_N)$ in $\bbA^n_X$. Clearly $W\cap V = U_Z^u$, and so $W \to X$ is
smooth along points of   $U_Z^u$. There exists a neighborhood $U_X^u
\subset W$ of $u$ containing $U_Z^u$ which is smooth over $X$. We
can take
   \[
   U_Z' = \bigcup_u U_Z^u, \quad U_X' = \bigcup_u U_X^u.\qedhere
   \]
\end{proof}

\begin{lemma}
Let $i_j:Z \to X_j$ be closed embeddings of algebraic stacks,
$j=1,2$.  There exist schemes $U_j$ and smooth surjective morphisms
$U_j \to
X_j$ such that
$$ Z\times_{X_1}U_1 \cong Z\times_{X_2}U_2. $$
\end{lemma}

\emph{Proof.} Let $V_j \to X_j$ be smooth and surjective. We have a pullback diagram
$$\xymatrix{ i_j^* V_j \ar[r]\ar[d]& V_j\ar[d]\\
                     Z\ar^{i_j}[r] & X_j.
}$$
Consider the fibered product $\tilde V = i_1^* V_1 \times_Z i_2^* V_2$. Applying the previous lemma with   $Z\subset X$ replaced by the closed embedding $i_1^*V_1 \subset V_1$, and  $U_Z \to Z$ replaced  by $\tilde V \to i_1^*V_1$, there is a Zariski-open covering $ \tilde V_1 \to \tilde V$ and a smooth $\tilde U_1\to X_1$ with $i_1^*\tilde U_1 \cong  \tilde V_1$. Applying the same procedure for  $i_1^*V_1 \subset V_1$, we obtain a  Zariski-open covering $ \tilde V_2 \to \tilde V$ and a smooth $\tilde U_2\to X_2$ with $i_2^*\tilde U_2 \cong  \tilde V_2$. Replacing $\tilde V_1$ and $\tilde V_2$ by a common Zariski-refinement $\tilde V_{12}$, and replacing $\tilde U_j$ by a suitable Zariski-refinements $U_j$ lifting  $\tilde V_{12}$, the lemma is proven. \qed

\subsection{The construction}\hfill

\begin{proof}[Proof of \ref{Prop:gluing}]
By the previous lemma, there is a choice of schemes $U_i$
and smooth
surjective morphisms $U_i \to X_i$ with $i_1^*U_1 \cong i_2^*U_2= U_Z$. Set
$R_i = U_i \times_{X_i} U_i$, so that $R_i \double U_i$ is a
presentation of $X_i$, and $i_1^*R_1 \cong i_2^*R_2= R_Z \double U_Z$
is a presentation of $Z$. Set $U = U_1 \cup_{U_Z} U_2$ and $R =
R_1\cup_{R_Z} R_2$. We have a groupoid $R \double U$, with a diagram
of groupoids
\begin{equation}\label{cubediagram}
\xymatrix{
                     &     R_Z\ar[ld]\doubledown\ar[rd] &\\
R_1\doubledown\ar[rd]& U_Z\ar[ld]\ar[rd]  & R_2\doubledown\ar[ld]\\
U_1\ar[rd]           & R \doubledown    & U_2\ar[ld] \\
                     & U                &}
                     \end{equation}
Let $X = [R\double U]$ be the quotient. We claim this is the desired stack.

Let $\cM$ be an algebraic stack. By definition, we have
$$\Hom(X,\cM)\ \ =\ \
\varprojlim\left(\cM(U)\ \double \ \cM(R) \right).$$

By Lemma \ref{Lem:gluing-is-coproduct}, we have

$$\cM(U)\ \ =\ \
\cM(U_1)\ \ \mathop{\times}\limits_{\cM(U_Z)}\ \ \cM(U_2)
\ \ = \ \  \varprojlim\left(\cM(U_1\sqcup U_2)\ \double \ \cM(U_Z) \right)
$$
and
$$\cM(R)\ \ =\ \
\cM(R_1)\ \ \mathop{\times}\limits_{\cM(R_Z)}\ \ \cM(R_2)
\ \ = \ \  \varprojlim\left(\cM(R_1\sqcup R_2)\ \double \ \cM(R_Z) \right).$$

Taking $\cM$-valued points in diagram (\ref{cubediagram}), we get
\begin{equation}\label{reversediagram}
\xymatrix{
                     &     \cM(R_Z)&\\
\cM(R_1)\ar[ur]& \cM(U_Z) \doubleup  & \cM(R_2)\ar[ul]\\
\cM(U_1)\doubleup\ar[ur]        &  & \cM(U_2)\doubleup\ar[ul] }
                     \end{equation}

We thus have

\begin{align*}
&\Hom(X,\cM)  = \\
 &= \varprojlim
\Big( \varprojlim\left(\cM(U_1\sqcup U_2)\double\cM(U_Z) \right)  \double
\varprojlim
\left(\cM(R_1\sqcup R_2)\double\cM(R_Z) \right)\Big) \\
 &= \varprojlim\, (\text{ diagram (\ref{reversediagram})\ })\\
 &= \varprojlim\left(
 \cM(U_1) \double \cM(R_1)
  \right)
 \mathop{\times}\limits_{\varprojlim
\left(\cM(U_Z) \rightrightarrows \cM(R_Z)\right)}
\left(\varprojlim
 \cM(U_2) \double \cM(R_2)
  \right)\\
&= \Hom(X_1,\cM) \ \mathop{\times}\limits_{\Hom(Z,\cM)}\ \Hom(X_2,\cM).
\qedhere\end{align*}
\end{proof}

\begin{proof}[Proof of corollary \ref{Cor:limits}]
Consider the morphisms
$$j_1 := i_1\sqcup i_2: Z \sqcup Z \to Y \ \ \text{and} \ \  j_2 :=
i_2\sqcup i_1: Z \sqcup Z \to Y.$$ These morphisms are closed
embeddings by the empty intersection hypothesis. Let $$\tilde X = Y
\mathop\cup\limits_{Z \sqcup Z}Y.$$ There is a canonical free $G =
\ZZ/2\ZZ$-action on $\tilde X$ arising from the action on the gluing
diagram. We claim that $$X = \left[\tilde X \ \big/ \  G \right]$$ is
the desired colimit. Indeed, $X$ is the colimit of the following
diagram:
  \[ 
 \xymatrix{
G \times (Z \sqcup Z) \doubleright \doubledown & Z \sqcup Z \doubledown \\
G \times (Y \sqcup Y) \doubleright                      & Y \sqcup Y.
}
 \]\end{proof}

\section{Taking roots of line bundles}

The following constrution is due independently to the authors and to C.~Cadman (\cite{Cadman1}, Section~2).

\subsection{The root of a line bundle}
If $L$ is a line bundle on a scheme $S$ and $d$ is a positive
integer, we will denote by $\radice{d}{L/S}$ the stack over $S$ of
$d\th$ roots of $L$; an object of $\radice{d}{L/S}$ over $T \to S$
is a line bundle $M$ over $T$, together with an isomorphism of
$M^{\otimes d}$ with the pullback $L_T$ of $L$ to $T$. The arrows
are defined in the obvious way. This stack $\radice{d}{L/S}$ is a
gerbe over $S$ banded by $\mu_d$; its cohomology class in the flat
cohomology group $\H^2(S, \mu_d)$ is obtained from the class $[L]
\in \H^1(S, \gm)$ via the boundary homomorphism $\partial\colon
\H^1(S, \gm) \to \H^2(S, \mu_d)$ obtained from the Kummer exact
sequence
   \[
   0 \longrightarrow \mu_{d,S}
   \longrightarrow\mathbb{G}_{\mathrm{m},S}
   \xrightarrow{(-)^d} \mathbb{G}_{\mathrm{m},S}
   \longrightarrow 0.
   \]

It is clear that if $T \to S$ is a morphism of schemes, then
$\radice{d}{L_T/T} = \radice{d}{L/S} \times_S T$.
The stack $\radice{d}{L/S}$ can be described directly as the quotient stack
$[L^0/\gm]$. Here $L^0$ is the total space of the $\gm$-bundle associated
with $L$ (or, in more down to earth terms, $L$ minus the zero section),
and the action of $\gm$ over $L^0$ is defined by the formula $\rho(\alpha)
x = \alpha^{-d}x$.
 Equivalently,  $$\radice{d}{L/S}\ =\ S \times_{\cB\GG_m}{\cB\GG_m},$$ where the fibered product is taken with respect to the classifying morphism $[L]: S \to {\cB\GG_m}$ on the left and the $d$-th power map ${\cB\GG_m}\to {\cB\GG_m}$ on the right. There is a universal line bundle $\mathcal{L}$ over
$\radice{d}{L/S}$, which is the quotient $[\mathbb{A}^1 \times L^0/ \gm]$
by the action defined by $\alpha(u, x) = (\alpha u, \rho(\alpha)x)$.

\subsection{The root of a line bundle with a section}\label{roots}

Now, given a section $\sigma\colon S \to L$ we can define a variant of
this, which we denote by $\radice{d}{(L, \sigma)/S}$, in which the objects
over a scheme $T \to S$ consist of triples
$$(M,\phi,\tau)$$
where
\begin{enumerate}
\item $M$ is a line bundle over $T$,
\item $\phi:M^{\otimes d} \simeq L_T$ is an isomorphism, and
\item $\tau$ is a section of $M$ such that $\phi(\tau^m) = \sigma$.
\end{enumerate}

 If $Y$ is the
scheme-theoretic zero locus of $\sigma$, then the restriction of the
stack $\radice{d}{(L, \sigma)/S}$ to $S \setminus Y$ is equal to $S
\setminus Y$. Its restriction to $Y$ is more interesting; it does
not coincide with $\radice{d}{L_{Y}/Y}$, because an object of the stack
$\radice{d}{(L, \sigma)/S}\mid_Y$ over $T \to Y$ consists of a
$d\th$ root $M$ of $L_T$, plus a section of $M$ whose $d\th$ power
is 0; but the section itself is not necessarily 0. However, the
morphism $\radice{d}{L_{Y}/Y} \to \radice{d}{(L, \sigma)/S}\mid_Y$
defined by sending a $d\th$ root $M$ of $L_T$ on a scheme $T$ over
$Y$ to the same $M$ together with the zero section is a closed
embedding, defined by a nilpotent sheaf of ideals on $\radice{d}{(L,
\sigma)/S}\mid_Y$. Thus $\radice{d}{(L, \sigma)/S}\mid_Y$ contains a
canonical gerbe banded by $\mu_{d}$, supported over the zero scheme
of $\sigma$. The forgetful map $\radice{d}{(L, \sigma)/S}\mid_Y\to
\radice{d}{L_{Y}/Y}$ identifies $\radice{d}{(L, \sigma)/S}\mid_Y$ as
the $d$-th infinitesimal neighborhood of $\radice{d}{L_{Y}/Y}$ in
its universal line bundle.

The stack $\radice{d}{(L, \sigma)/S}$ can also be decribed as a quotient
stack. Consider the universal line bundle $\mathcal{L} = [\mathbb{A}^1
\times L^0/ \gm]$ described above, and the morphism $\Phi\colon \mathcal{L}
\to L$ induced by the $\gm$-invariant morphism $\mathbb{A}^1 \times L^0 \to
L$ defined by $(u, x) \mapsto u^d x$; then
   \[
   \radice{d}{(L, \sigma)/S} =
   \Phi^{-1} \sigma(S).
   \]
In other words: if we call
$V_\sigma\subseteq
\mathbb{A}^1 \times L^0$ the inverse image in $\mathbb{A}^1 \times L^0$ of
the embedding $\sigma \colon S \into L$; then
   \[
   \radice{d}{(L, \sigma)/S} = [V_\sigma/\gm].
   \]

In particular, assume that $L = \cO$, so that $\sigma$ is a
regular function on $S$. In this case $L^0 = S \times \gm$, and if we
denote by $W_\sigma = \mathbb{A}^1 \times S$ the subscheme defined by the
equation $t^d - \sigma(x) = 0$, where $t$ is a coordinate on
$\mathbb{A}^1$, then there is an isomorphism $W_\sigma \times \gm \simeq
V_\sigma$.

An equivalent description exists here too: Let $\cU  = [\AA^1/\GG_m]$, the classifying stack for line bundles with section. Then  $$\radice{d}{(L, \sigma)/S} = S \times_\cU \cU,$$ where the map on the left is the classifying map and the map $\cU \to \cU$ on the right is the $d$-th power map.

\section{Rigidification}\label{Sec:rigidification}
\subsection{The setup}
We recall the concept of \emph{rigidification} of an algebraic
stack, as presented in \cite{ACV}, see also the related treatment in
\cite{Romagny}.

Let $H$ be a flat finitely presented separated group scheme over a
base  scheme $\bbS$, ${\cX}$ an algebraic  stack over $\bbS$. We say that $\cX$ has \emph{an $H$-2-structure}
if for each object $\xi \in {\cX}(T)$ there is an embedding
$$
\iota_\xi\colon H(T)\into \Aut_T (\xi),
$$
which is compatible with
pullback, in the following sense: given two objects $\xi\in
{\cX}(T)$ and $\eta\in{\cX}(T)$, and an arrow $\phi\colon \xi \to
\eta$ in ${\cX}$ over a morphism of schemes $f\colon S \to T$, the
natural pullback homomorphisms
$$\phi^*\colon \Aut_T (\eta)\to \Aut_S (\xi)$$ and $$f^* \colon H(T) \to H(S)$$
commute with the embedding, that is, $\iota_\xi f^* = \phi^*\iota_\eta$.

This condition can also be expressed as follows. Let $\phi\colon
\xi \to \eta$ be an arrow in ${\cX}$ over a morphism of schemes
$f\colon S \to
T$, and
$g \in H(T)$. Then
the diagram
   \[
   \xymatrix{
   \xi\ar[r]^\phi\ar[d]^{f^*g} & \eta\ar[d]^g \\
   \xi\ar[r]^\phi & \eta
   }
   \]
commutes. In particular, by taking $\xi = \eta$ and $\phi$ to be in
$\Aut_S(\xi)$, we see that $H(S)$ must be in the center of
$\Aut_S(\xi)$; in particular, $H$ is an abelian group
scheme. 

The simplest example of such a situation is when $\cX \to T$ is a
gerbe banded by $H$; in this case the embedding $H(S)\into \Aut_S
(\xi)$ is an isomorphism of group schemes.

Then we have the following result.

\begin{theorem}(\cite{ACV}, Theorem~5.1.5)
There is a smooth surjective finitely presented morphism of algebraic stacks $\cX \to \cX\thickslash H$ satisfying the following properties:

\begin{enumerate}
\item For any object $\xi\in \cX(T)$ with image $\eta\in \cX\thickslash H(T)$, we have that $H(T)$ lies in the kernel of $\Aut_T( \xi) \to \Aut_T(\eta)$.
\item The morphism $\cX \to \cX\thickslash H$ is universal for  morphisms of stacks $\cX \to \cY$ satisfying (1) above.
\item If $T$ is the spectrum of an algebraically closed field, then in (1) above,  we have
   \[
   \Aut_T(\eta) = \Aut_T( \xi)/H(T).
   \]
\item A moduli space for $\cX$ is also a moduli space for $\cX\thickslash H$.
\end{enumerate}

Furthermore, if $\cX$ is a Deligne--Mumford stack, then $\cX\thickslash H$ is
also a Deligne--Mumford stack and the morphism $\cX \to \cX\thickslash H$ is
\'etale.

\end{theorem}

The notation in \cite{ACV} is $\cX^{H}$; here we adopt the better
notation $\cX\thickslash H$ proposed by Romagny in \cite{Romagny}.
This stack $\cX\thickslash H$ is called the \emph{$H$-rigidification
of $\cX$}. For example, if $\cG \to T$ is a gerbe banded by $H$,
then $\cG\thickslash H$ is isomorphic to $T$.

The stack $\cX\thickslash H$ is obtained as the fppf stackification
of a prestack $\cX^{H}_{\mathrm{pre}}$, that has the same objects as
$\cX$; this has the property that for any object $\xi$ of $\cX$ over
an $\bbS$-scheme $T$, the sheaf of automorphisms of
$\cAut_{T,\cX^{H}_{\mathrm{pre}}}(\xi)$ in $\cX^{H}_{\mathrm{pre}}$
is the quotient sheaf of $\cAut_{T,\cX}(\xi)$ by the normal subgroup
sheaf $H_{T}$.

\subsection{Moduli interpretation of the stack $\cX\thickslash H$}\label{subsection-moduli-interpretation}

The construction of rigidification is functorial, in the
sense described below. First let $\phi\colon  \cX \to \cY$ be a morphism of algebraic
stacks endowed with $H$-2-structures. We say that $\phi$ is \emph{$H$-2-equivariant} if for each
$\bbS$-scheme $T$ and object $\xi$ of $\cX(T)$, the homomorphism of
group-schemes
   \[
   H(T) \into \Aut_{T,\cX}(\xi) \longrightarrow \Aut_{T,\cY}(\phi(\xi))
   \]
defined by $\phi$ coincides with the given embedding $H(T)
\hookrightarrow \Aut_{T,\cY}(\phi\xi)$.

Define a 2-category $\cX\slashsecond H$ over the category of schemes over $\bbS$, as follows.

\begin{enumerate}

\item An object over a scheme $T$ is a pair $(\cG,\phi)$, where  $\cG \to T$ is a gerbe banded by $H$, and $\cG \stackrel{\phi}{\to} \cX$ is an $H$-2-equivariant morphism of fibered categories.

\item A morphism  $(F,\rho): (\cG,\phi) \to (\cG',\phi')$ consists of a morphism $F : \cG \to \cG'$ over some $f: T \to T'$, compatible with the bands, and a 2-morphism $\rho: \phi\to  \phi'\circ F$ making the following diagram commutative:
$$\xymatrix{ \cG \ar^{F}[rr]\ar[rd]_\phi && \cG' \ar[ld]^{\phi'} \\
& \cX
}$$

\item A 2-arrow $(F,\rho) \to  (F_{1},\rho_{1})$ is a usual  2-arrow $\sigma: F \to F_{1}$ compatible with $\rho$ and $\rho_{1}$ in the sense that the following diagram is commutative:
   \[
   \xymatrix{ & \phi \ar[dl]_{{\rho}} \ar[dr]^{\rho_{1}}\\
    \phi'\circ F \ar[rr]^{\phi'(\sigma)} && \phi'\circ F_{1}
   }
   \]
\end{enumerate}

From Lemma~\ref{lem:representable->1-category} we see that  $\cX\slashsecond H$ is equivalent to a 1-category, that we denote by $\cX\slashprime H$. This is easily checked to be a category fibered in groupoids over the category of schemes over $\bbS$.

\begin{proposition}\label{prop:equivalence-rigidification}
There is an equivalence of fibered categories between $\cX\thickslash H$ and $\cX\slashprime H$.
\end{proposition}

\begin{proof}
Given an $H$-2-equivariant morphism $\cG \to \cX$, where $\cG \to T$ is a gerbe banded by $H$, there is an induced morphism $T = \cG\thickslash H \to \cX\thickslash H$. This gives a function from the objects of $\cX\slashprime H$ to the objects of $\cX\thickslash H$, that extends to a functor in the obvious way.

In the other direction, given an object $\xi$ of $\cX\thickslash H(T)$, consider the fibered product $\cG := T \times_{\cX\thickslash H} \cX \to T$. If $U$ is a scheme over $T$, an object of $\cG(U)$ is a pair $(\zeta,\alpha)$, where $\zeta$ is an object of $\cX(U)$, and $\alpha$ is an isomorphism between $\zeta$ and the pullback $\xi_{U}$ in $\cX\thickslash H(U)$. We claim that $\cG$ is a gerbe over $T$, and there is a unique banding of $\cG$ by $H$ making the projection $\cG \to \cX$ $H$-2-equivariant.

Both statements are fppf local on $T$, so we may assume that the given object $\xi$ of $\cX\thickslash H (T)$ comes from an object $\widetilde{\xi}$ of $\cX(T)$. Then the pair $(\widetilde{\xi}_{U},\mathrm{id})$ is an object of $\cG(U)$, showing the existence of local sections. Two objects $(\zeta_{1}, \alpha_{1})$ and $(\zeta_{2}, \alpha_{2})$ in $\cG(U)$ are locally isomorphic, because the morphism of sheaves of sets $\cHom_{U,\cX}(\zeta_{1}, \zeta_{2}) \to \cHom_{U,\cX\thickslash H}(\zeta_{1}, \zeta_{2})$ is an $H$-torsor.

Also, given an object $(\zeta,\alpha)$ of $\cG(U)$, its automorphism group in $\cG(U)$ is the kernel of the homomorphism $\Aut_{U,\cX}(\zeta) \to \Aut_{U,\cX\thickslash H}(\zeta)$, that is exactly $H(U) \subseteq \Aut_{U,\cX}(\zeta)$.

So $\cG$ is a gerbe banded by $H_{T}$. This function from the objects of $\cX$ to the objects of $\cX\slashprime H$ extends to a functor in the obvious way.

It is immediate to see that the composition $\cX\thickslash H \to \cX\slashprime H \to \cX\thickslash H$ is isomorphic to the identity. Let us show that the composition $\cX\slashprime H \to \cX\thickslash H \to \cX\slashprime H$ is also isomorphic to the identity. Given a gerbe $\cG \to T$ and an $H$-2-equivariant morphism $\cG \to \cX$, the induced morphism $\cG \to T \times_{\cX\thickslash H} \cX$ is a morphism of gerbes banded by $H$; and any such morphism is an isomorphism.
\end{proof}

\subsection{The rigidification as a quotient} Here is another way to think of the rigidification. Given the collection of data $\iota_\xi : H(T) \hookrightarrow \Aut_T(\xi)$ for all objects $\xi$ of $\cX$ required for forming the rigidification, we have an associated action $$\cB H \times \cX \to \cX,$$ defined as follows. First, given a scheme $T$, an object of $(\cB H \times \cX)\  (T)$ consists of  a principal $H$-bundle $P \to T$ together with an object $\xi\in \cX(T)$. We need to form a new object $P\star \xi \in \cX(T)$. We do this as follows: the pullback $\xi_P$ of $\xi$ to $P$ admits a left diagonal action of $H$ coming from the two actions on $P$ (inverted, as to make it into a left action) and on $\xi$. By the descent axiom for $\cX$ there is a quotient object on $T$ which we call  $P\star \xi \in \cX(T)$. The assumptions on  $\iota_\xi$ guarantee that the formation of $P\star \xi$ is functorial, giving the required morphism $\cB H \times \cX \to \cX$. The particular case $\cX = \cB H$ shows that $\cB H$ is indeed a group stack, and in general one shows the morphism $\cB H \times \cX \to \cX$ is an action. The morphism $\cX \to \cX \thickslash H$ is easily seen to be invariant. From the fact that $\iota_\xi$ is injective one obtaines that  $$\cB H \times \cX\ \  \lrar\ \  \cX\mathop\times\limits_{\cX \thickslash H}\cX $$
is an isomorphism, so this action is free and, whatever the 2-categorical quotient  should mean, up to equivalence we obtain
$$ \cX / \cB H \ \ \simeq \ \ \cX \thickslash H.$$

This is certainly in agreement with our previous moduli interpretation of $\cX \thickslash H$: a principal $\cB H$-bundle $\cG \to T$ is simply a gerbe banded by $H$, and saying that $\cG \to \cX$ is equivariant translates to our requirement on the map of automorphisms to coincide with the homomorphisms $\iota_\xi$.

The case of non-central actions is more complicated and is worked
out in the appendix to \cite{AOV}.

\end{document}